%Paper: alg-geom/9603020
%From: hosono@sci1.sci.toyama-u.ac.jp
%Date: Mon, 25 Mar 1996 01:32:31 +0900
%Date (revised): Mon, 1 Apr 1996 11:57:51 +0900

\input harvmac

\def\frac#1#2{{#1\over#2}}

\def\journal#1&#2(#3){\unskip, #1~\bf #2 \rm(19#3) }
\def\andjournal#1&#2(#3){\sl #1~\bf #2 \rm (19#3) }

\catcode`\@=11\def\slash#1{\mathord{\mathpalette\c@ncel{#1}}}
\overfullrule=0pt
\def\steepslash{\c@ncel}
\def\frac#1#2{{#1\over #2}}

\def\:{\!:\!}
\def\inbar{\,\vrule height1.5ex width.4pt depth0pt}
\def\IQ{\relax\,\hbox{$\inbar\kern-.3em{\rm Q}$}}
\def\IB{\relax{\rm I\kern-.18em B}}
\def\IC{\relax\hbox{$\inbar\kern-.3em{\rm C}$}}
\def\IP{\relax{\rm I\kern-.18em P}}
\def\IR{\relax{\rm I\kern-.18em R}}
\def\ZZ{\relax\ifmmode\mathchoice
{\hbox{Z\kern-.4em Z}}{\hbox{Z\kern-.4em Z}}
{\lower.9pt\hbox{Z\kern-.4em Z}}
{\lower1.2pt\hbox{Z\kern-.4em Z}}\else{Z\kern-.4em Z}\fi}

\catcode`\@=12

%                      Zeitschriften:
\def\npb#1(#2)#3{{ Nucl. Phys. }{B#1} (#2) #3}
\def\plb#1(#2)#3{{ Phys. Lett. }{#1B} (#2) #3}
\def\pla#1(#2)#3{{ Phys. Lett. }{#1A} (#2) #3}
\def\prl#1(#2)#3{{ Phys. Rev. Lett. }{#1} (#2) #3}
\def\mpla#1(#2)#3{{ Mod. Phys. Lett. }{A#1} (#2) #3}
\def\ijmpa#1(#2)#3{{ Int. J. Mod. Phys. }{A#1} (#2) #3}
\def\cmp#1(#2)#3{{ Comm. Math. Phys. }{#1} (#2) #3}
\def\cqg#1(#2)#3{{ Class. Quantum Grav. }{#1} (#2) #3}
\def\jmp#1(#2)#3{{ J. Math. Phys. }{#1} (#2) #3}
\def\anp#1(#2)#3{{ Ann. Phys. }{#1} (#2) #3}
\def\prd#1(#2)#3{{ Phys. Rev. } {D{#1}} (#2) #3}
\def\ptp#1(#2)#3{{ Progr. Theor. Phys. }{#1} (#2) #3}
\def\aom#1(#2)#3{{ Ann. Math. }{#1} (#2) #3}

\def\C{{\bf C}}

\def\Z{{\bf Z}}
\def\P{{\bf P}}

\def\Z{{\bf Z}}

\def\cicy#1(#2|#3)#4{\left(\matrix{#2}\right|\!\!
                     \left|\matrix{#3}\right)^{{#4}}_{#1}}

\def\ra{\rightarrow}

\def\Box{{\,\lower0.9pt\vbox{\hrule
\hbox{\vrule height 0.2 cm \hskip 0.2 cm
\vrule height 0.2 cm}\hrule}\,}}

\global\newcount\thmno \global\thmno=0
\def\definition#1{\global\advance\thmno by1
\bigskip\noindent{\bf Definition \secsym\the\thmno. }{\it #1}
\par\nobreak\medskip\nobreak}
\def\question#1{\global\advance\thmno by1
\bigskip\noindent{\bf Question \secsym\the\thmno. }{\it #1}
\par\nobreak\medskip\nobreak}
\def\theorem#1{\global\advance\thmno by1
\bigskip\noindent{\bf Theorem \secsym\the\thmno. }{\it #1}
\par\nobreak\medskip\nobreak}
\def\proposition#1{\global\advance\thmno by1
\bigskip\noindent{\bf Proposition \secsym\the\thmno. }{\it #1}
\par\nobreak\medskip\nobreak}
\def\corollary#1{\global\advance\thmno by1
\bigskip\noindent{\bf Corollary \secsym\the\thmno. }{\it #1}
\par\nobreak\medskip\nobreak}
\def\lemma#1{\global\advance\thmno by1
\bigskip\noindent{\bf Lemma \secsym\the\thmno. }{\it #1}
\par\nobreak\medskip\nobreak}
\def\conjecture#1{\global\advance\thmno by1
\bigskip\noindent{\bf Conjecture \secsym\the\thmno. }{\it #1}
\par\nobreak\medskip\nobreak}
\def\exercise#1{\global\advance\thmno by1
\bigskip\noindent{\bf Exercise \secsym\the\thmno. }{\it #1}
\par\nobreak\medskip\nobreak}
\def\remark#1{\global\advance\thmno by1
\bigskip\noindent{\bf Remark \secsym\the\thmno. }{\it #1}
\par\nobreak\medskip\nobreak}
\def\problem#1{\global\advance\thmno by1
\bigskip\noindent{\bf Problem \secsym\the\thmno. }{\it #1}
\par\nobreak\medskip\nobreak}
\def\others#1#2{\global\advance\thmno by1
\bigskip\noindent{\bf #1 \secsym\the\thmno. }{\it #2}
\par\nobreak\medskip\nobreak}

\def\thmlab#1{\xdef
#1{\secsym\the\thmno}\writedef{#1\leftbracket#1}\wrlabeL{#1=#1}}
%
% redefine \newsec so that all \thmno set to zero in a new section
%
\catcode`\@=11
\def\s@csym{}
\def\newsec#1{\global\advance\secno by1%
{\toks0{#1}\message{(\the\secno. \the\toks0)}}%
%\ifx\answ\bigans \vfill\eject \else \bigbreak\bigskip \fi  %if desired
\global\subsecno=0\thmno=0%
\eqnres@t\let\s@csym\secsym\xdef\secn@m{\the\secno}\noindent
{\bf\hyperdef\hypernoname{section}{\the\secno}{\the\secno.} #1}%
\writetoca{{\string\hyperref{}{section}{\the\secno}{\the\secno.}} {#1}}%
\par\nobreak\medskip\nobreak}
\catcode`\@=12
\def\fivepoint{\def\rm{\fam0\fiverm}% switch to fivepoint font
\textfont0=\fiverm \scriptfont0=\fiverm \scriptscriptfont0=\fiverm
\textfont1=\fivei \scriptfont1=\fivei \scriptscriptfont1=\fivei
\textfont2=\fivesy \scriptfont2=\fivesy \scriptscriptfont2=\fivesy
\textfont\itfam=\fivei \def\it{\fam\itfam\fiveit}\def\sl{\fam\slfam\fivesl}%
\textfont\bffam=\fivebf \def\bf{\fam\bffam\fivebf}\rm}

\lref\AspinwallGross{P. Aspinwall and M. Gross, {\sl
Heterotic-Heterotic String Duality and Multiple K3 Fibrations},
hep-th/9602118.}
\lref\KlemmSchimmrigk{A. Klemm and R. Schimmrigk,
{\sl Landau-Ginzburg String Vacua}, CERN-TH-6459/92, Nucl. Phys. B.}
\lref\KreuzerSkarke{M. Kreuzer and H. Skarke, Nucl. Phys. B388 (1993) 113}
\lref\KachruVafa{ S. Kachru and C. Vafa,
       {\sl Exact Results for N=2 Compactifications of Heterotic Strings},
      hep-th/9505105.}
\lref\KachruSilverstein{S. Kachru and E. Silverstein,
{\sl N=1 Dual String Pairs and Gaugino Condensation},
hep-th/9511228.}
\lref\KlemmLercheMayr{A. Klemm, W. Lerche and P. Mayr,
       {\sl K3-Fibrations and Heterotic-Type II String Duality},
	     hep-th/9506112.}
\lref\Yonemura{ T. Yonemura,
T\^ohoku Math. J. {\bf 42}(1990),351.}
\lref\HKTYI{S.Hosono, A.Klemm, S.Theisen and S.-T.Yau,
	    Commun. Math. Phys. 167 (1995) 301.}
\lref\CFKM{P.Candelas,  A.Font, S.Katz and D.Morrison,
		  Nucl.Phys.{\bf B416}(1994)481.}
\lref\HosonoLianYau{ S. Hosono, B. Lian and  S.T. Yau, {\sl
GKZ-Generalized Hypergeometric
Systems in Mirror Symmetry of Calabi-Yau Hypersurfaces}, Harvard Univ.
preprint, alg-geom/9511001, to appear in CMP 1996.}
\lref\Oguiso{K. Oguiso, Int. J. Math. 4 (1993) 439.}

\Title{alg-geom/9603020}
{\vbox{
 \centerline{Calabi-Yau Varieties and Pencils of K3 Surfaces}
    } }

   \centerline{S. Hosono$^1$\footnote{$^\dagger$}
   {email: hosono@sci.toyama-u.ac.jp},
   B.H. Lian$^2$\footnote{$^\ddagger$}{email: lian@max.math.brandeis.edu}
   and S.-T. Yau$^3$\footnote{$^\diamond$}{email: yau@math.harvard.edu}}

   \bigskip\centerline{
   \vbox{
   \hbox{ $^1$ Department of Mathematics}
   \hbox{ \hskip30pt Toyama University}
   \hbox{ \hskip28pt Toyama 930, Japan} }
   \hskip0.5cm
   \vbox{
   \hbox{ $^2$ Department of Mathematics}
   \hbox{ \hskip30pt Brandeis University}
   \hbox{ \hskip27pt Waltham, MA 02154} } }
   \bigskip
   \centerline{
   \vbox{
   \hbox{ $^3$ Department of Mathematics}
   \hbox{ \hskip30pt Harvard University}
   \hbox{ \hskip23pt Cambridge, MA 02138} } }

\vskip.3in

Abstract: In this note, we give a list of Calabi-Yau hypersurfaces
in weighted projective 4-spaces with the property that
a hypersurface contains naturally a pencil of K3 variety.
For completeness we also obtain a similar list in the case
K3 hypersurfaces in weighted projective 3-spaces.
The first list significantly enlarges the list of
K3-fibrations of \KlemmLercheMayr~ which has been obtained
on some assumptions on the weights.
Our lists are expected to correspond to examples of the so-called
heterotic-type II duality \KachruVafa\KachruSilverstein.

\Date{3/2/96} %If you remove this line, you lose page numbers !!!!!

\newsec{Problems}

Let $w_1,..,w_{n+1}$ be positive integers, and put
$d:=\sum w_i$. We call the
weight vector $\hat{w}=(w_1,..,w_{n+1})$ admissible if
the generic weighted degree $d$
hypersurface in $\C^{n+1}$ is smooth away from the origin. This means that
the weighted projectivized hypersurface in $\P[w]$ is transversal, ie. it only
acquires
singularities from the ambient space $\P[\hat w]$.
For $n=3$, there is a list of admissible weights of
Reid-Yonemura (see \Yonemura). For $n=4$, there is a list
admissible weights obtained
by Klemm-Schimmrigk \KlemmSchimmrigk~ and Kreuzer-Skarke \KreuzerSkarke.

Given an admissible weight $\hat{w}=(w_1,..,w_{n+1})$
 we can consider in $\P[\hat w]$ the generic
Calabi-Yau variety given by
\eqn\dumb{\hat X_a=\{z|\sum_{\hat{w}\cdot\nu=d} a_\nu z^\nu=0\}.}
Suppose we intersect this variety with the coordinate hyperplane
$z_{n+1}=0$.

\problem{When is $X_a:=\hat{X}_a\cap\{z_{n+1}=0\}$ isomorphic to
a transversal Calabi-Yau variety?}

Note that by permuting the weights, this includes the cases
$\hat X_a\cap\{z_i=0\}$ for any $i$.
More generally,

\problem{When is there a 1-parameter family of hypersurfaces
$Z_\lambda$ such that $X_{a,\lambda}:=\hat X_a\cap Z_\lambda$
is isomorphic to a transversal Calabi-Yau variety?}

For $n=4$ and with some assumptions on the weights, a short list
of such cases has been tabulated in \KlemmLercheMayr.
We say that $\nu$ is compatible with the weight $\hat w$ if
$\hat w\cdot\nu=d$.

 Let $w=(w_1,..,w_n)$ and
$\bar{w}$ its normalization, ie. $\bar{w}:=(w_1/\delta_1,..,w_n/\delta_n)$
where $\delta_i:=lcm(\rho_1,..,\hat{\rho_i},..,\rho_n)$ and
 $\rho_i:=gcd(w_1,..,\hat{w_i},..,w_n)$.
It is well known that $\phi:\P[w]\ra\P[\bar w]$
is an isomorphism under the normalization map
$(z_1,..,z_n)\mapsto(z_1^{\rho_1},..,z_n^{\rho_n})$.
It is easy to show that $\delta_1\rho_1=\cdots=\delta_n\rho_n$;
we call this integer $k$.

We require that the image $\bar{X}_a=\phi X_a$ is a transversal
Calabi-Yau variety in $\P[\bar w]$. If $x_1,..,x_n$ are the
quasi-homogeneous coordinates of $\P[\bar w]$, then a Calabi-Yau variety
can be written as
\eqn\dumb{\bar X_b=\{x|\sum_{\bar w\cdot\mu=\bar d} b_\mu x^\mu=0\},}
where $\bar d:=\sum_{i=1}^n \bar w_i=\sum w_i/\delta_i$.
Pulling this back by the normalization map,
we see that
\eqn\dumb{
\phi^{-1}\bar X_b=
\{z|\sum_{\bar w\cdot\mu=\bar d} b_\mu \prod z_i^{\rho_i\mu_i}=0\}\subset
\P[w].}
If we require that $\phi^{-1}\bar X_b=X_a$ for some $a$,
then we conclude that
\item{(a)} every $\mu\in\Z_+^n$ compatible
with $\bar w$ has $\sum w_i\rho_i\mu_i=d$, and
$\nu=(\rho_1\mu_1,..,\rho_n\mu_n,0)$ is an exponent compatible with
$\hat w$.
\item{(b)} every exponent
$\nu\in\Z_+^{n+1}$ compatible with $\hat w$ having $\nu_{n+1}=0$
is of the form $\nu=(\rho_1\mu_1,..,\rho_n\mu_n,0)$ for some
exponent $\mu$ compatible with $\bar w$.

We claim that $d=k\bar d$. This follows from (a) and the fact
that $k=\rho_i\delta_i$.
Since $\mu=(1,..,1)$ is compatible with $\bar w$, it
 follows from (a) that $w\cdot\rho=d$.
Thus our task is to search through the list of
normalized admissible weights
$\hat w$ ($n=4$) satisfying
\eqn\dumb{\eqalign{
(i)&~~\hat w\cdot\nu=d,~~\nu_5=0~\Rightarrow\rho_i|\nu_i~~\forall i\cr
(ii)&~~ w\cdot\rho=d\cr
(iii)&~~\bar w~~admissible}}
On the last condition, we will check that $\bar w$ be in the
Reid-Yonemura list. It is also clear that (i)--(iii) implies that
$\hat X_a\cap\{z_{n+1}\}$ is isomorphic to $\bar X_a$ in the
admissible $\P[\bar w]$. Our computer search shows that there are
628 admissible weights $\hat w$ of length 5 satisfying (i)--(iii).

{\it Example:} Take $\hat w=(42,27,8,4,3)$, $d=84$. We consider the
intersection $X:=\hat X\cap\{z_3=0\}$. Then $\rho=(1,1,3,1)$,
and so condition (ii) holds.
The normalized weight of $w=(42,27,4,3)$ becomes $\bar w=(14,9,4,1)$,
which is an admissible weight of length 4 (see \Yonemura),
and so condition (iii) holds.
The equations for $X$ in $\P[\hat w]$ is $z_3=0$ plus that of $\hat X$.
The latter is given by the generic sum
of the monomials with admissible exponent $\nu$ with $\nu_5=0$. There
are exactly 24 such exponents:
\eqn\dumb{\eqalign{
& z_5^{28},~ z_5^{24}z_4^3,~ z_5^{20}z_4^6,~ z_5^{16}z_4^9,~
z_5^{12} z_4^{12},~ z_5^8z_4^{15},~ z_5^4z_4^{18},~ z_4^{21},~ z_5^{19}z_2, \cr
& z_5^{15}z_4^3z_2,~ z_5^{11}z_4^6z_2,~ z_5^7z_4^9z_2,~
z_5^3z_4^{12}z_2,~ z_5^{10}z_2^2,~ z_5^6z_4^3z_2^2,~\cr
& z_5^2z_4^6z_2^2,~ z_5z_2^3,~ z_5^{14}z_1,~ z_5^{10}z_4^3z_1,~
z_5^6z_4^6z_1,~ z_5^2z_4^9z_1,~ z_5^5z_2z_1,~
z_5z_4^3z_2z_1,~ z_1^2.}}
Condition (i) holds because the exponent $\nu_4$ of $z_4$ is always
a multiple of $\rho_3=3$.
The equation for the isomorphic image $\bar X$ of $X$ in $\P[14,9,4,1]$
is the generic sum of the above monomials with the replacement,
$z_1\mapsto x_1$, $z_2\mapsto x_2$, $z_4\mapsto x_3^3$,
$z_5\mapsto x_4$.

We note that given an admissible weight $\hat w$, the Calabi-Yau
varieties in $\P[\hat w]$ can give two distinct transversal
Calabi-Yau varieties when intersect with two different coordinate
hyperplanes $z_i=0$.

\subsec{the second problem}

We consider our second problem under the following assumption.
We assume that $Z_\lambda$ is of the form $\lambda_1z_{n+1}=\lambda_2 p(z)$
where $\lambda=[\lambda_1,\lambda_2]$ is regarded as a point
in $\P^1$, and $p(z)$ a fixed nonzero
quasi-homogeneous polynomial
independent of $z_{n+1}$ and has degree $w_{n+1}$. When
$\lambda_2=0$ this reduces to the case in the first problem.
This generalization turns out to require just some minor modification.
Specifically, in addition to conditions (i)--(iii), we must
require that the weight component
\item{(iv)} $w_{n+1}$ can be
partitioned by the components $w_1,..,w_n$.

This is true iff $p$ exists.
Note that as $\lambda$ varies the intersections
$\hat X_a\cap Z_\lambda$ form a pencil of codimension one subvarieties
in $\hat X_a$.
In the case of $n=4$ we require that
they are transversal K3 varieties when $\lambda_1\neq0$.
In our list of 628 cases above,  we find that all of them admit this
description hence enlarging the list of \KlemmLercheMayr.

The table given in the appendix is
the list of the 628 cases. The number
denoted $i$ between 1 and 5 in the table indicates $Z_\lambda$ is
of the form $\lambda_1 z_i=\lambda_2 p(z)$ as in the case $i=n+1$
discussed above. Some of the examples in this list have been
studied in great details in the context of mirror symmetry
(see for example \HKTYI\CFKM\HosonoLianYau), and in connection with
string duality in \KachruVafa\AspinwallGross~ and others.

We note that the conditions we impose
in our method for enumerating K3 pencils are only sufficient
but not necessary. There is in fact a criterion given in
\Oguiso~ for K3 pencils using the intersection ring of
the Calabi-Yau variety. In fact in \HosonoLianYau~ (see the conclusion
section there) we have already
used this criterion to
give a few examples of K3 pencils in which we have computed
the intersection ring.
For example, the Calabi-Yau hypersurfaces in $\P[8,3,3,1,1]$
was found to have a K3 pencil according to the criterion of \Oguiso,
but this example fails to satisfy conditions (i)-(iii) above.
In \HosonoLianYau, we have also given an algorithm for
computing the intersection ring of Calabi-Yau hypersurfaces
in weighted projective spaces. This algorithm can in principle
be carried out for all of the list \KlemmSchimmrigk, and be used
to check the criterion above. But the actual computation can
be enormous.

For completeness, we also do the case of $n=3$. Thus we search
through the list of transversal K3 hypersurfaces in \Yonemura~
which admits a pencil of elliptic curves in one of the
following transversal weighted projective spaces
$\P[1,1,1], \P[2,1,1], \P[3,2,1]$.
The $n=3$ analogues of conditions (i)-(iii)
are satisfied by
18 admissible weights, and all
of them satisfy condition (iv).

\vskip.3in
{\it Acknowledgement:} We thank S. Kachru and C. Vafa for helpful
discussions, and C. Doran for pointing out an error in the early
version of this paper.

$${\fivepoint{
\vbox{\offinterlineskip\tabskip=0pt
\halign{\strut
\vrule#
&~~$#$~~\hfil
&~~$#$~~\hfil
&~~$#$~~\hfil
&\vrule#
&~~$#$~~\hfil
&~~$#$~~\hfil
&~~$#$~~\hfil
&\vrule#\cr
\noalign{\hrule}
 &\hat w=(w_1,..,w_4)& i&~~\bar w~~ &
 &\hat w=(w_1,..,w_4)& i&~~\bar w~~ &\cr
\noalign{\hrule}
& (4, 3, 3, 2) & 1 & (2, 1, 1) & & (4, 3, 3, 2) & 2 & (3, 2, 1) &\cr
& (2, 2, 1, 1) & 3 & (1, 1, 1) & & (4, 4, 3, 1) & 3 & (1, 1, 1) &\cr
& (4, 2, 1, 1) & 3 & (2, 1, 1) & & (6, 3, 2, 1) & 3 & (2, 1, 1) &\cr
& (10, 5, 4, 1) & 3 & (2, 1, 1) & & (6, 4, 1, 1) & 3 & (3, 2, 1) &\cr
& (9, 6, 2, 1) & 3 & (3, 2, 1) & & (12, 8, 3, 1) & 3 & (3, 2, 1) &\cr
& (21, 14, 6, 1) & 3 & (3, 2, 1) & & (3, 3, 2, 1) & 3 & (1, 1, 1) &\cr
& (9, 4, 3, 2) & 2 & (3, 2, 1) & & (8, 4, 3, 1) & 3 & (2, 1, 1) &\cr
& (12, 7, 3, 2) & 2 & (2, 1, 1) & & (18, 11, 4, 3) & 2 & (3, 2, 1) &\cr
& (15, 10, 4, 1) & 3 & (3, 2, 1) & & (18, 12, 5, 1) & 3 & (3, 2, 1) &\cr
\noalign{\hrule} }}}}$$

\listrefs

\newsec{Appendix}

$${\fivepoint{
\vbox{\offinterlineskip\tabskip=0pt
\halign{\strut\vrule#
&~~$#$~~\hfil
&~~$#$~~\hfil
&~~$#$~~\hfil
&~~$#$~~\hfil
&~~$#$~~\hfil
&\vrule#\cr
\noalign{\hrule}
&Euler~\#& h^{1,1}& \hat w=(w_1,..,w_5) &i &~~~\bar w&\cr
\noalign{\hrule}
& 480 & 287 & ( 882, 588, 251, 36, 7 ) & 3 & ( 21, 14, 6, 1 ) &\cr
& 376 & 201 & ( 280, 140, 109, 16, 15 ) & 3 & ( 14, 7, 4, 3 ) &\cr
& 324 & 212 & ( 630, 420, 179, 24, 7 ) & 3 & ( 15, 10, 4, 1 ) &\cr
& 256 & 147 & ( 200, 100, 77, 15, 8 ) & 3 & ( 10, 5, 3, 2 ) &\cr
& 240 & 173 & ( 504, 336, 143, 18, 7 ) & 3 & ( 12, 8, 3, 1 ) &\cr
& 216 & 141 & ( 200, 100, 79, 16, 5 ) & 3 & ( 10, 5, 4, 1 ) &\cr
& 192 & 110 & ( 60, 60, 43, 9, 8 ) & 3 & ( 5, 5, 3, 2 ) &\cr
& 180 & 114 & ( 143, 110, 44, 30, 3 ) & 4 & ( 13, 10, 4, 3 ) &\cr
& 180 & 114 & ( 130, 100, 40, 27, 3 ) & 4 & ( 13, 10, 4, 3 ) &\cr
& 168 & 95 & ( 144, 67, 48, 20, 9 ) & 2 & ( 12, 5, 4, 3 ) &\cr
\noalign{\hrule}
& 160 & 115 & ( 160, 80, 61, 15, 4 ) & 3 & ( 8, 4, 3, 1 ) &\cr
& 160 & 115 & ( 160, 80, 63, 12, 5 ) & 3 & ( 8, 4, 3, 1 ) &\cr
& 156 & 86 & ( 77, 56, 42, 30, 5 ) & 4 & ( 11, 8, 6, 5 ) &\cr
& 156 & 86 & ( 66, 48, 36, 25, 5 ) & 4 & ( 11, 8, 6, 5 ) &\cr
& 144 & 131 & ( 378, 252, 107, 12, 7 ) & 3 & ( 9, 6, 2, 1 ) &\cr
& 144 & 98 & ( 162, 99, 32, 27, 4 ) & 3 & ( 18, 11, 4, 3 ) &\cr
& 144 & 91 & ( 80, 56, 32, 21, 3 ) & 4 & ( 10, 7, 4, 3 ) &\cr
& 120 & 86 & ( 48, 48, 35, 9, 4 ) & 3 & ( 4, 4, 3, 1 ) &\cr
& 120 & 69 & ( 100, 35, 32, 25, 8 ) & 3 & ( 20, 8, 7, 5 ) &\cr
& 120 & 69 & ( 54, 42, 25, 24, 5 ) & 3 & ( 9, 7, 5, 4 ) &\cr
\noalign{\hrule}
& 120 & 65 & ( 60, 40, 36, 35, 9 ) & 3 & ( 12, 9, 8, 7 ) &\cr
& 120 & 65 & ( 48, 32, 28, 27, 9 ) & 4 & ( 12, 9, 8, 7 ) &\cr
& 112 & 76 & ( 98, 49, 24, 21, 4 ) & 3 & ( 14, 7, 4, 3 ) &\cr
& 112 & 63 & ( 55, 30, 28, 20, 7 ) & 3 & ( 11, 7, 6, 4 ) &\cr
& 112 & 63 & ( 44, 24, 21, 16, 7 ) & 3 & ( 11, 7, 6, 4 ) &\cr
& 108 & 60 & ( 50, 30, 25, 24, 21 ) & 1 & ( 25, 10, 8, 7 ) &\cr
& 108 & 60 & ( 25, 25, 20, 16, 14 ) & 1 & ( 25, 10, 8, 7 ) &\cr
& 96 & 167 & ( 225, 200, 150, 24, 1 ) & 4 & ( 9, 8, 6, 1 ) &\cr
& 96 & 167 & ( 216, 192, 144, 23, 1 ) & 4 & ( 9, 8, 6, 1 ) &\cr
& 96 & 87 & ( 120, 60, 47, 8, 5 ) & 3 & ( 6, 3, 2, 1 ) &\cr
\noalign{\hrule} }}}}$$
\vfill\eject

$${\fivepoint{
\vbox{\offinterlineskip\tabskip=0pt
\halign{\strut\vrule#
&~~$#$~~\hfil
&~~$#$~~\hfil
&~~$#$~~\hfil
&~~$#$~~\hfil
&~~$#$~~\hfil
&\vrule#\cr
\noalign{\hrule}
&Euler~\#& h^{1,1}& \hat w=(w_1,..,w_5) &i &~~~\bar w&\cr
\noalign{\hrule}
& 96 & 79 & ( 88, 64, 21, 16, 3 ) & 3 & ( 11, 8, 3, 2 ) &\cr
& 96 & 65 & ( 63, 42, 35, 24, 4 ) & 4 & ( 9, 6, 5, 4 ) &\cr
& 96 & 59 & ( 99, 44, 22, 18, 15 ) & 2 & ( 33, 22, 6, 5 ) &\cr
& 96 & 59 & ( 44, 33, 33, 12, 10 ) & 2 & ( 33, 22, 6, 5 ) &\cr
& 96 & 59 & ( 44, 39, 22, 15, 12 ) & 1 & ( 22, 13, 5, 4 ) &\cr
& 96 & 59 & ( 56, 33, 20, 12, 11 ) & 2 & ( 14, 11, 5, 3 ) &\cr
& 96 & 59 & ( 42, 22, 15, 11, 9 ) & 2 & ( 14, 11, 5, 3 ) &\cr
& 96 & 57 & ( 44, 32, 24, 15, 5 ) & 4 & ( 11, 8, 6, 5 ) &\cr
& 96 & 55 & ( 38, 24, 19, 18, 15 ) & 1 & ( 19, 8, 6, 5 ) &\cr
& 96 & 55 & ( 19, 19, 16, 12, 10 ) & 1 & ( 19, 8, 6, 5 ) &\cr
\noalign{\hrule}
& 84 & 104 & ( 294, 196, 56, 39, 3 ) & 4 & ( 21, 14, 4, 3 ) &\cr
& 84 & 54 & ( 36, 27, 27, 10, 8 ) & 2 & ( 27, 18, 5, 4 ) &\cr
& 84 & 54 & ( 38, 33, 19, 15, 9 ) & 1 & ( 19, 11, 5, 3 ) &\cr
& 84 & 54 & ( 22, 19, 19, 10, 6 ) & 2 & ( 19, 11, 5, 3 ) &\cr
& 84 & 50 & ( 36, 31, 18, 15, 8 ) & 2 & ( 6, 5, 4, 3 ) &\cr
& 84 & 50 & ( 34, 21, 18, 17, 12 ) & 1 & ( 17, 7, 6, 4 ) &\cr
& 84 & 50 & ( 17, 17, 14, 12, 8 ) & 1 & ( 17, 7, 6, 4 ) &\cr
& 80 & 68 & ( 112, 56, 32, 21, 3 ) & 4 & ( 14, 7, 4, 3 ) &\cr
& 80 & 51 & ( 32, 16, 15, 12, 5 ) & 3 & ( 8, 5, 4, 3 ) &\cr
& 72 & 68 & ( 108, 49, 36, 20, 3 ) & 2 & ( 9, 5, 3, 1 ) &\cr
\noalign{\hrule}
& 72 & 68 & ( 108, 53, 36, 15, 4 ) & 2 & ( 9, 5, 3, 1 ) &\cr
& 72 & 65 & ( 36, 36, 25, 8, 3 ) & 3 & ( 3, 3, 2, 1 ) &\cr
& 72 & 59 & ( 56, 35, 18, 14, 3 ) & 3 & ( 8, 5, 3, 2 ) &\cr
& 72 & 57 & ( 50, 35, 20, 12, 3 ) & 4 & ( 10, 7, 4, 3 ) &\cr
& 72 & 50 & ( 72, 32, 16, 15, 9 ) & 2 & ( 24, 16, 5, 3 ) &\cr
& 72 & 50 & ( 44, 27, 20, 9, 8 ) & 2 & ( 11, 9, 5, 2 ) &\cr
& 72 & 49 & ( 34, 30, 17, 12, 9 ) & 1 & ( 17, 10, 4, 3 ) &\cr
& 72 & 49 & ( 28, 24, 15, 8, 5 ) & 3 & ( 7, 6, 5, 2 ) &\cr
& 72 & 49 & ( 28, 11, 11, 10, 6 ) & 2 & ( 14, 11, 5, 3 ) &\cr
& 72 & 49 & ( 20, 17, 17, 8, 6 ) & 2 & ( 17, 10, 4, 3 ) &\cr
\noalign{\hrule} }}}}$$
\vfill\eject

$${\fivepoint{
\vbox{\offinterlineskip\tabskip=0pt
\halign{\strut\vrule#
&~~$#$~~\hfil
&~~$#$~~\hfil
&~~$#$~~\hfil
&~~$#$~~\hfil
&~~$#$~~\hfil
&\vrule#\cr
\noalign{\hrule}
&Euler~\#& h^{1,1}& \hat w=(w_1,..,w_5) &i &~~~\bar w&\cr
\noalign{\hrule}
& 72 & 48 & ( 33, 24, 18, 10, 5 ) & 4 & ( 11, 8, 6, 5 ) &\cr
& 72 & 47 & ( 60, 21, 16, 15, 8 ) & 3 & ( 20, 8, 7, 5 ) &\cr
& 72 & 47 & ( 27, 21, 12, 10, 5 ) & 4 & ( 9, 7, 5, 4 ) &\cr
& 72 & 46 & ( 33, 18, 14, 12, 7 ) & 3 & ( 11, 7, 6, 4 ) &\cr
& 72 & 44 & ( 24, 16, 14, 9, 9 ) & 4 & ( 12, 9, 8, 7 ) &\cr
& 72 & 44 & ( 24, 21, 20, 12, 7 ) & 2 & ( 7, 6, 5, 3 ) &\cr
& 72 & 44 & ( 18, 15, 14, 9, 7 ) & 3 & ( 7, 6, 5, 3 ) &\cr
& 72 & 44 & ( 21, 18, 16, 9, 8 ) & 3 & ( 8, 7, 6, 3 ) &\cr
& 64 & 47 & ( 36, 16, 15, 8, 5 ) & 3 & ( 9, 5, 4, 2 ) &\cr
& 64 & 43 & ( 22, 12, 8, 7, 7 ) & 4 & ( 11, 7, 6, 4 ) &\cr
\noalign{\hrule}
& 60 & 194 & ( 465, 248, 186, 30, 1 ) & 4 & ( 15, 8, 6, 1 ) &\cr
& 60 & 194 & ( 450, 240, 180, 29, 1 ) & 4 & ( 15, 8, 6, 1 ) &\cr
& 60 & 59 & ( 90, 55, 16, 15, 4 ) & 3 & ( 18, 11, 4, 3 ) &\cr
& 60 & 49 & ( 25, 25, 12, 10, 3 ) & 3 & ( 5, 5, 3, 2 ) &\cr
& 60 & 44 & ( 22, 16, 12, 5, 5 ) & 4 & ( 11, 8, 6, 5 ) &\cr
& 54 & 56 & ( 35, 35, 21, 12, 2 ) & 4 & ( 5, 5, 3, 2 ) &\cr
& 50 & 44 & ( 35, 25, 20, 12, 3 ) & 4 & ( 7, 5, 4, 3 ) &\cr
& 48 & 83 & ( 156, 91, 39, 24, 2 ) & 4 & ( 12, 7, 3, 2 ) &\cr
& 48 & 59 & ( 96, 40, 32, 21, 3 ) & 4 & ( 12, 5, 4, 3 ) &\cr
& 48 & 53 & ( 52, 40, 16, 9, 3 ) & 4 & ( 13, 10, 4, 3 ) &\cr
\noalign{\hrule}
& 48 & 43 & ( 28, 21, 21, 10, 4 ) & 2 & ( 21, 14, 5, 2 ) &\cr
& 48 & 43 & ( 32, 27, 16, 15, 6 ) & 1 & ( 16, 9, 5, 2 ) &\cr
& 48 & 41 & ( 36, 21, 12, 8, 7 ) & 2 & ( 9, 7, 3, 2 ) &\cr
& 48 & 41 & ( 27, 14, 9, 7, 6 ) & 2 & ( 9, 7, 3, 2 ) &\cr
& 48 & 41 & ( 22, 10, 9, 9, 4 ) & 3 & ( 11, 9, 5, 2 ) &\cr
& 48 & 39 & ( 18, 14, 8, 5, 5 ) & 4 & ( 9, 7, 5, 4 ) &\cr
& 48 & 39 & ( 26, 18, 15, 13, 6 ) & 1 & ( 13, 6, 5, 2 ) &\cr
& 48 & 39 & ( 40, 24, 21, 20, 15 ) & 1 & ( 20, 8, 7, 5 ) &\cr
& 48 & 39 & ( 32, 21, 16, 15, 12 ) & 1 & ( 16, 7, 5, 4 ) &\cr
& 48 & 39 & ( 15, 15, 14, 12, 4 ) & 1 & ( 15, 7, 6, 2 ) &\cr
\noalign{\hrule} }}}}$$
\vfill\eject

$${\fivepoint{
\vbox{\offinterlineskip\tabskip=0pt
\halign{\strut\vrule#
&~~$#$~~\hfil
&~~$#$~~\hfil
&~~$#$~~\hfil
&~~$#$~~\hfil
&~~$#$~~\hfil
&\vrule#\cr
\noalign{\hrule}
&Euler~\#& h^{1,1}& \hat w=(w_1,..,w_5) &i &~~~\bar w&\cr
\noalign{\hrule}
& 48 & 39 & ( 13, 13, 12, 10, 4 ) & 1 & ( 13, 6, 5, 2 ) &\cr
& 48 & 38 & ( 40, 21, 16, 12, 7 ) & 2 & ( 10, 7, 4, 3 ) &\cr
& 48 & 37 & ( 24, 12, 10, 9, 5 ) & 3 & ( 8, 5, 4, 3 ) &\cr
& 48 & 36 & ( 36, 24, 16, 15, 5 ) & 4 & ( 9, 6, 5, 4 ) &\cr
& 48 & 35 & ( 16, 15, 12, 12, 5 ) & 2 & ( 5, 4, 3, 3 ) &\cr
& 48 & 35 & ( 12, 10, 7, 7, 6 ) & 3 & ( 7, 6, 5, 3 ) &\cr
& 48 & 35 & ( 12, 10, 9, 9, 5 ) & 2 & ( 5, 4, 3, 3 ) &\cr
& 44 & 51 & ( 70, 35, 18, 14, 3 ) & 3 & ( 10, 5, 3, 2 ) &\cr
& 42 & 55 & ( 49, 35, 21, 12, 2 ) & 4 & ( 7, 5, 3, 2 ) &\cr
& 40 & 69 & ( 110, 55, 33, 20, 2 ) & 4 & ( 10, 5, 3, 2 ) &\cr
\noalign{\hrule}
& 40 & 41 & ( 40, 25, 20, 12, 3 ) & 4 & ( 8, 5, 4, 3 ) &\cr
& 36 & 116 & ( 133, 114, 76, 18, 1 ) & 4 & ( 7, 6, 4, 1 ) &\cr
& 36 & 116 & ( 126, 108, 72, 17, 1 ) & 4 & ( 7, 6, 4, 1 ) &\cr
& 36 & 38 & ( 28, 21, 21, 8, 6 ) & 2 & ( 21, 14, 4, 3 ) &\cr
& 36 & 38 & ( 26, 24, 13, 9, 6 ) & 1 & ( 13, 8, 3, 2 ) &\cr
& 36 & 38 & ( 16, 13, 13, 6, 4 ) & 2 & ( 13, 8, 3, 2 ) &\cr
& 36 & 35 & ( 21, 18, 10, 6, 5 ) & 3 & ( 7, 6, 5, 2 ) &\cr
& 36 & 34 & ( 66, 31, 15, 12, 8 ) & 2 & ( 11, 5, 4, 2 ) &\cr
& 36 & 34 & ( 22, 15, 12, 11, 6 ) & 1 & ( 11, 5, 4, 2 ) &\cr
& 36 & 34 & ( 11, 11, 10, 8, 4 ) & 1 & ( 11, 5, 4, 2 ) &\cr
\noalign{\hrule}
& 32 & 103 & ( 102, 85, 68, 16, 1 ) & 4 & ( 6, 5, 4, 1 ) &\cr
& 32 & 103 & ( 96, 80, 64, 15, 1 ) & 4 & ( 6, 5, 4, 1 ) &\cr
& 32 & 46 & ( 70, 35, 20, 12, 3 ) & 4 & ( 14, 7, 4, 3 ) &\cr
& 32 & 33 & ( 16, 8, 6, 5, 5 ) & 4 & ( 8, 5, 4, 3 ) &\cr
& 24 & 80 & ( 216, 144, 43, 27, 2 ) & 3 & ( 12, 8, 3, 1 ) &\cr
& 24 & 49 & ( 100, 40, 33, 25, 2 ) & 3 & ( 10, 5, 4, 1 ) &\cr
& 24 & 48 & ( 55, 40, 12, 10, 3 ) & 3 & ( 11, 8, 3, 2 ) &\cr
& 24 & 47 & ( 72, 35, 24, 9, 4 ) & 2 & ( 6, 3, 2, 1 ) &\cr
& 24 & 41 & ( 105, 42, 30, 28, 5 ) & 3 & ( 15, 6, 5, 4 ) &\cr
& 24 & 41 & ( 90, 36, 25, 24, 5 ) & 3 & ( 15, 6, 5, 4 ) &\cr
\noalign{\hrule} }}}}$$
\vfill\eject

$${\fivepoint{
\vbox{\offinterlineskip\tabskip=0pt
\halign{\strut\vrule#
&~~$#$~~\hfil
&~~$#$~~\hfil
&~~$#$~~\hfil
&~~$#$~~\hfil
&~~$#$~~\hfil
&\vrule#\cr
\noalign{\hrule}
&Euler~\#& h^{1,1}& \hat w=(w_1,..,w_5) &i &~~~\bar w&\cr
\noalign{\hrule}
& 24 & 38 & ( 42, 21, 9, 8, 4 ) & 4 & ( 14, 7, 4, 3 ) &\cr
& 24 & 38 & ( 63, 28, 15, 14, 6 ) & 2 & ( 21, 14, 5, 2 ) &\cr
& 24 & 38 & ( 32, 20, 9, 8, 3 ) & 3 & ( 8, 5, 3, 2 ) &\cr
& 24 & 38 & ( 27, 18, 16, 8, 3 ) & 3 & ( 9, 8, 6, 1 ) &\cr
& 24 & 38 & ( 32, 21, 20, 7, 4 ) & 2 & ( 8, 7, 5, 1 ) &\cr
& 24 & 38 & ( 24, 15, 14, 7, 3 ) & 3 & ( 8, 7, 5, 1 ) &\cr
& 24 & 38 & ( 33, 16, 9, 8, 6 ) & 2 & ( 11, 8, 3, 2 ) &\cr
& 24 & 36 & ( 35, 20, 12, 10, 3 ) & 3 & ( 7, 4, 3, 2 ) &\cr
& 24 & 34 & ( 27, 18, 15, 8, 4 ) & 4 & ( 9, 6, 5, 4 ) &\cr
& 24 & 33 & ( 18, 7, 7, 6, 4 ) & 2 & ( 9, 7, 3, 2 ) &\cr
\noalign{\hrule}
& 24 & 33 & ( 25, 12, 10, 10, 3 ) & 2 & ( 5, 3, 2, 2 ) &\cr
& 24 & 33 & ( 27, 12, 10, 6, 5 ) & 3 & ( 9, 5, 4, 2 ) &\cr
& 24 & 33 & ( 20, 9, 8, 8, 3 ) & 2 & ( 5, 3, 2, 2 ) &\cr
& 24 & 32 & ( 14, 12, 5, 5, 4 ) & 3 & ( 7, 6, 5, 2 ) &\cr
& 24 & 32 & ( 32, 15, 12, 8, 5 ) & 2 & ( 8, 5, 3, 2 ) &\cr
& 24 & 32 & ( 30, 17, 12, 9, 4 ) & 2 & ( 5, 3, 2, 2 ) &\cr
& 24 & 31 & ( 28, 15, 12, 8, 5 ) & 2 & ( 7, 5, 3, 2 ) &\cr
& 24 & 30 & ( 54, 25, 12, 9, 8 ) & 2 & ( 9, 4, 3, 2 ) &\cr
& 24 & 29 & ( 15, 15, 12, 10, 8 ) & 1 & ( 15, 6, 5, 4 ) &\cr
& 24 & 28 & ( 24, 18, 17, 9, 4 ) & 3 & ( 4, 3, 3, 2 ) &\cr
\noalign{\hrule}
& 24 & 27 & ( 27, 18, 12, 10, 5 ) & 4 & ( 9, 6, 5, 4 ) &\cr
& 24 & 27 & ( 30, 14, 12, 9, 7 ) & 2 & ( 10, 7, 4, 3 ) &\cr
& 24 & 27 & ( 8, 6, 6, 5, 5 ) & 4 & ( 5, 4, 3, 3 ) &\cr
& 20 & 50 & ( 91, 56, 18, 14, 3 ) & 3 & ( 13, 8, 3, 2 ) &\cr
& 18 & 53 & ( 63, 49, 21, 12, 2 ) & 4 & ( 9, 7, 3, 2 ) &\cr
& 16 & 31 & ( 18, 8, 5, 5, 4 ) & 3 & ( 9, 5, 4, 2 ) &\cr
& 16 & 29 & ( 20, 16, 9, 8, 3 ) & 3 & ( 5, 4, 3, 2 ) &\cr
& 12 & 41 & ( 54, 33, 9, 8, 4 ) & 4 & ( 18, 11, 4, 3 ) &\cr
& 12 & 36 & ( 26, 21, 15, 13, 3 ) & 1 & ( 13, 7, 5, 1 ) &\cr
& 12 & 36 & ( 16, 15, 15, 12, 2 ) & 2 & ( 15, 8, 6, 1 ) &\cr
\noalign{\hrule} }}}}$$
\vfill\eject

$${\fivepoint{
\vbox{\offinterlineskip\tabskip=0pt
\halign{\strut\vrule#
&~~$#$~~\hfil
&~~$#$~~\hfil
&~~$#$~~\hfil
&~~$#$~~\hfil
&~~$#$~~\hfil
&\vrule#\cr
\noalign{\hrule}
&Euler~\#& h^{1,1}& \hat w=(w_1,..,w_5) &i &~~~\bar w&\cr
\noalign{\hrule}
& 12 & 36 & ( 14, 13, 13, 10, 2 ) & 2 & ( 13, 7, 5, 1 ) &\cr
& 8 & 29 & ( 28, 16, 9, 8, 3 ) & 3 & ( 7, 4, 3, 2 ) &\cr
& 6 & 23 & ( 21, 10, 9, 6, 5 ) & 2 & ( 7, 5, 3, 2 ) &\cr
& 0 & 251 & ( 903, 602, 258, 42, 1 ) & 4 & ( 21, 14, 6, 1 ) &\cr
& 0 & 251 & ( 882, 588, 252, 41, 1 ) & 4 & ( 21, 14, 6, 1 ) &\cr
& 0 & 131 & ( 253, 138, 92, 22, 1 ) & 4 & ( 11, 6, 4, 1 ) &\cr
& 0 & 131 & ( 242, 132, 88, 21, 1 ) & 4 & ( 11, 6, 4, 1 ) &\cr
& 0 & 121 & ( 153, 136, 102, 16, 1 ) & 4 & ( 9, 8, 6, 1 ) &\cr
& 0 & 119 & ( 210, 105, 84, 20, 1 ) & 4 & ( 10, 5, 4, 1 ) &\cr
& 0 & 119 & ( 200, 100, 80, 19, 1 ) & 4 & ( 10, 5, 4, 1 ) &\cr
\noalign{\hrule}
& 0 & 89 & ( 225, 150, 45, 28, 2 ) & 4 & ( 15, 10, 3, 2 ) &\cr
& 0 & 89 & ( 96, 80, 48, 15, 1 ) & 4 & ( 6, 5, 3, 1 ) &\cr
& 0 & 89 & ( 90, 75, 45, 14, 1 ) & 4 & ( 6, 5, 3, 1 ) &\cr
& 0 & 83 & ( 252, 168, 71, 7, 6 ) & 3 & ( 6, 4, 1, 1 ) &\cr
& 0 & 77 & ( 70, 56, 42, 13, 1 ) & 4 & ( 5, 4, 3, 1 ) &\cr
& 0 & 77 & ( 65, 52, 39, 12, 1 ) & 4 & ( 5, 4, 3, 1 ) &\cr
& 0 & 71 & ( 52, 52, 39, 12, 1 ) & 4 & ( 4, 4, 3, 1 ) &\cr
& 0 & 71 & ( 48, 48, 36, 11, 1 ) & 4 & ( 4, 4, 3, 1 ) &\cr
& 0 & 65 & ( 168, 112, 32, 21, 3 ) & 4 & ( 21, 14, 4, 3 ) &\cr
& 0 & 59 & ( 165, 110, 30, 22, 3 ) & 3 & ( 15, 10, 3, 2 ) &\cr
\noalign{\hrule}
& 0 & 59 & ( 150, 100, 27, 20, 3 ) & 3 & ( 15, 10, 3, 2 ) &\cr
& 0 & 55 & ( 80, 40, 31, 5, 4 ) & 3 & ( 4, 2, 1, 1 ) &\cr
& 0 & 55 & ( 147, 98, 36, 7, 6 ) & 3 & ( 21, 14, 6, 1 ) &\cr
& 0 & 55 & ( 98, 63, 24, 7, 4 ) & 3 & ( 14, 9, 4, 1 ) &\cr
& 0 & 55 & ( 80, 56, 21, 8, 3 ) & 3 & ( 10, 7, 3, 1 ) &\cr
& 0 & 55 & ( 70, 49, 18, 7, 3 ) & 3 & ( 10, 7, 3, 1 ) &\cr
& 0 & 41 & ( 24, 24, 17, 4, 3 ) & 3 & ( 2, 2, 1, 1 ) &\cr
& 0 & 39 & ( 50, 25, 16, 5, 4 ) & 3 & ( 10, 5, 4, 1 ) &\cr
& 0 & 39 & ( 75, 35, 24, 10, 6 ) & 3 & ( 15, 7, 6, 2 ) &\cr
& 0 & 39 & ( 35, 20, 12, 5, 3 ) & 3 & ( 7, 4, 3, 1 ) &\cr
\noalign{\hrule} }}}}$$
\vfill\eject

$${\fivepoint{
\vbox{\offinterlineskip\tabskip=0pt
\halign{\strut\vrule#
&~~$#$~~\hfil
&~~$#$~~\hfil
&~~$#$~~\hfil
&~~$#$~~\hfil
&~~$#$~~\hfil
&\vrule#\cr
\noalign{\hrule}
&Euler~\#& h^{1,1}& \hat w=(w_1,..,w_5) &i &~~~\bar w&\cr
\noalign{\hrule}
& 0 & 39 & ( 45, 25, 16, 10, 4 ) & 3 & ( 9, 5, 4, 2 ) &\cr
& 0 & 38 & ( 60, 25, 20, 12, 3 ) & 4 & ( 12, 5, 4, 3 ) &\cr
& 0 & 35 & ( 63, 28, 18, 14, 3 ) & 2 & ( 21, 14, 6, 1 ) &\cr
& 0 & 35 & ( 63, 28, 18, 14, 3 ) & 3 & ( 9, 4, 3, 2 ) &\cr
& 0 & 35 & ( 28, 21, 21, 12, 2 ) & 2 & ( 21, 14, 6, 1 ) &\cr
& 0 & 35 & ( 28, 21, 21, 12, 2 ) & 4 & ( 4, 3, 3, 2 ) &\cr
& 0 & 35 & ( 28, 27, 14, 12, 3 ) & 1 & ( 14, 9, 4, 1 ) &\cr
& 0 & 35 & ( 40, 21, 12, 7, 4 ) & 2 & ( 10, 7, 3, 1 ) &\cr
& 0 & 35 & ( 30, 14, 9, 7, 3 ) & 2 & ( 10, 7, 3, 1 ) &\cr
& 0 & 34 & ( 20, 14, 8, 3, 3 ) & 4 & ( 10, 7, 4, 3 ) &\cr
\noalign{\hrule}
& 0 & 31 & ( 60, 24, 16, 15, 5 ) & 4 & ( 15, 6, 5, 4 ) &\cr
& 0 & 31 & ( 20, 20, 12, 5, 3 ) & 3 & ( 4, 4, 3, 1 ) &\cr
& 0 & 31 & ( 16, 16, 9, 4, 3 ) & 3 & ( 4, 4, 3, 1 ) &\cr
& 0 & 31 & ( 22, 18, 12, 11, 3 ) & 1 & ( 11, 6, 4, 1 ) &\cr
& 0 & 31 & ( 24, 19, 15, 12, 2 ) & 2 & ( 5, 4, 2, 1 ) &\cr
& 0 & 31 & ( 24, 23, 12, 10, 3 ) & 2 & ( 5, 4, 2, 1 ) &\cr
& 0 & 31 & ( 40, 35, 24, 15, 6 ) & 3 & ( 8, 7, 6, 3 ) &\cr
& 0 & 31 & ( 35, 16, 15, 10, 4 ) & 2 & ( 7, 4, 3, 2 ) &\cr
& 0 & 31 & ( 16, 10, 7, 7, 2 ) & 3 & ( 8, 7, 5, 1 ) &\cr
& 0 & 31 & ( 12, 11, 11, 8, 2 ) & 2 & ( 11, 6, 4, 1 ) &\cr
\noalign{\hrule}
& 0 & 29 & ( 45, 20, 10, 9, 6 ) & 2 & ( 15, 10, 3, 2 ) &\cr
& 0 & 29 & ( 24, 15, 12, 5, 4 ) & 2 & ( 6, 5, 3, 1 ) &\cr
& 0 & 29 & ( 18, 10, 9, 5, 3 ) & 2 & ( 6, 5, 3, 1 ) &\cr
& 0 & 29 & ( 13, 12, 12, 9, 2 ) & 1 & ( 3, 2, 2, 1 ) &\cr
& 0 & 27 & ( 20, 15, 15, 6, 4 ) & 2 & ( 15, 10, 3, 2 ) &\cr
& 0 & 23 & ( 28, 21, 14, 12, 9 ) & 1 & ( 14, 7, 4, 3 ) &\cr
& 0 & 23 & ( 14, 9, 7, 6, 6 ) & 1 & ( 7, 3, 2, 2 ) &\cr
& 0 & 23 & ( 9, 9, 8, 6, 4 ) & 1 & ( 9, 4, 3, 2 ) &\cr
& 0 & 23 & ( 9, 9, 8, 6, 4 ) & 3 & ( 4, 3, 3, 2 ) &\cr
& 0 & 23 & ( 7, 7, 6, 4, 4 ) & 1 & ( 7, 3, 2, 2 ) &\cr
\noalign{\hrule} }}}}$$
\vfill\eject

$${\fivepoint{
\vbox{\offinterlineskip\tabskip=0pt
\halign{\strut\vrule#
&~~$#$~~\hfil
&~~$#$~~\hfil
&~~$#$~~\hfil
&~~$#$~~\hfil
&~~$#$~~\hfil
&\vrule#\cr
\noalign{\hrule}
&Euler~\#& h^{1,1}& \hat w=(w_1,..,w_5) &i &~~~\bar w&\cr
\noalign{\hrule}
& 0 & 22 & ( 24, 10, 9, 6, 5 ) & 2 & ( 8, 5, 3, 2 ) &\cr
& 0 & 18 & ( 18, 12, 8, 5, 5 ) & 4 & ( 9, 6, 5, 4 ) &\cr
& 0 & 18 & ( 20, 8, 7, 7, 6 ) & 3 & ( 10, 7, 4, 3 ) &\cr
& -4 & 26 & ( 35, 12, 10, 10, 3 ) & 2 & ( 7, 3, 2, 2 ) &\cr
& -8 & 29 & ( 40, 20, 9, 8, 3 ) & 3 & ( 10, 5, 3, 2 ) &\cr
& -8 & 25 & ( 24, 16, 15, 5, 4 ) & 3 & ( 6, 5, 4, 1 ) &\cr
& -12 & 38 & ( 26, 20, 8, 3, 3 ) & 4 & ( 13, 10, 4, 3 ) &\cr
& -12 & 30 & ( 20, 15, 15, 8, 2 ) & 2 & ( 15, 10, 4, 1 ) &\cr
& -12 & 30 & ( 20, 15, 15, 8, 2 ) & 4 & ( 4, 3, 3, 2 ) &\cr
& -12 & 30 & ( 22, 21, 11, 9, 3 ) & 1 & ( 11, 7, 3, 1 ) &\cr
\noalign{\hrule}
& -12 & 30 & ( 14, 11, 11, 6, 2 ) & 2 & ( 11, 7, 3, 1 ) &\cr
& -12 & 25 & ( 10, 10, 4, 3, 3 ) & 4 & ( 5, 5, 3, 2 ) &\cr
& -12 & 24 & ( 18, 13, 12, 9, 2 ) & 2 & ( 3, 3, 2, 1 ) &\cr
& -12 & 16 & ( 10, 9, 6, 6, 5 ) & 1 & ( 5, 3, 2, 2 ) &\cr
& -16 & 23 & ( 70, 28, 20, 15, 7 ) & 2 & ( 14, 7, 4, 3 ) &\cr
& -16 & 23 & ( 56, 21, 16, 12, 7 ) & 2 & ( 14, 7, 4, 3 ) &\cr
& -20 & 15 & ( 14, 6, 5, 5, 4 ) & 3 & ( 7, 5, 3, 2 ) &\cr
& -24 & 110 & ( 144, 128, 96, 15, 1 ) & 4 & ( 9, 8, 6, 1 ) &\cr
& -24 & 77 & ( 72, 60, 48, 11, 1 ) & 4 & ( 6, 5, 4, 1 ) &\cr
& -24 & 60 & ( 40, 40, 30, 9, 1 ) & 4 & ( 4, 4, 3, 1 ) &\cr
\noalign{\hrule}
& -24 & 51 & ( 84, 49, 21, 12, 2 ) & 4 & ( 12, 7, 3, 2 ) &\cr
& -24 & 29 & ( 20, 7, 7, 6, 2 ) & 2 & ( 10, 7, 3, 1 ) &\cr
& -24 & 27 & ( 16, 10, 4, 3, 3 ) & 4 & ( 8, 5, 3, 2 ) &\cr
& -24 & 26 & ( 36, 16, 9, 8, 3 ) & 2 & ( 12, 8, 3, 1 ) &\cr
& -24 & 26 & ( 36, 16, 9, 8, 3 ) & 3 & ( 9, 4, 3, 2 ) &\cr
& -24 & 26 & ( 28, 15, 8, 5, 4 ) & 2 & ( 7, 5, 2, 1 ) &\cr
& -24 & 26 & ( 21, 10, 6, 5, 3 ) & 2 & ( 7, 5, 2, 1 ) &\cr
& -24 & 25 & ( 48, 23, 15, 6, 4 ) & 2 & ( 8, 5, 2, 1 ) &\cr
& -24 & 23 & ( 12, 6, 5, 5, 2 ) & 3 & ( 6, 5, 3, 1 ) &\cr
& -24 & 22 & ( 36, 17, 9, 6, 4 ) & 2 & ( 6, 3, 2, 1 ) &\cr
\noalign{\hrule} }}}}$$
\vfill\eject

$${\fivepoint{
\vbox{\offinterlineskip\tabskip=0pt
\halign{\strut\vrule#
&~~$#$~~\hfil
&~~$#$~~\hfil
&~~$#$~~\hfil
&~~$#$~~\hfil
&~~$#$~~\hfil
&\vrule#\cr
\noalign{\hrule}
&Euler~\#& h^{1,1}& \hat w=(w_1,..,w_5) &i &~~~\bar w&\cr
\noalign{\hrule}
& -24 & 22 & ( 21, 9, 8, 6, 4 ) & 3 & ( 7, 4, 3, 2 ) &\cr
& -24 & 21 & ( 40, 15, 12, 8, 5 ) & 2 & ( 10, 5, 3, 2 ) &\cr
& -24 & 21 & ( 18, 12, 11, 4, 3 ) & 3 & ( 3, 2, 2, 1 ) &\cr
& -24 & 20 & ( 42, 14, 12, 9, 7 ) & 2 & ( 14, 7, 4, 3 ) &\cr
& -24 & 20 & ( 10, 4, 4, 3, 3 ) & 4 & ( 5, 3, 2, 2 ) &\cr
& -24 & 18 & ( 18, 12, 10, 5, 3 ) & 3 & ( 6, 5, 4, 1 ) &\cr
& -24 & 17 & ( 20, 15, 10, 9, 6 ) & 1 & ( 10, 5, 3, 2 ) &\cr
& -24 & 15 & ( 16, 6, 5, 5, 4 ) & 3 & ( 8, 5, 3, 2 ) &\cr
& -24 & 12 & ( 6, 5, 5, 4, 4 ) & 2 & ( 5, 3, 2, 2 ) &\cr
& -28 & 17 & ( 14, 10, 8, 3, 3 ) & 4 & ( 7, 5, 4, 3 ) &\cr
\noalign{\hrule}
& -30 & 24 & ( 15, 15, 8, 5, 2 ) & 3 & ( 3, 3, 2, 1 ) &\cr
& -30 & 23 & ( 15, 15, 9, 4, 2 ) & 4 & ( 5, 5, 3, 2 ) &\cr
& -30 & 17 & ( 15, 9, 8, 4, 3 ) & 3 & ( 5, 4, 3, 1 ) &\cr
& -32 & 87 & ( 136, 68, 51, 16, 1 ) & 4 & ( 8, 4, 3, 1 ) &\cr
& -32 & 87 & ( 128, 64, 48, 15, 1 ) & 4 & ( 8, 4, 3, 1 ) &\cr
& -32 & 29 & ( 52, 32, 9, 8, 3 ) & 3 & ( 13, 8, 3, 2 ) &\cr
& -32 & 19 & ( 16, 9, 8, 4, 3 ) & 2 & ( 4, 3, 2, 1 ) &\cr
& -32 & 17 & ( 16, 10, 8, 3, 3 ) & 4 & ( 8, 5, 4, 3 ) &\cr
& -36 & 148 & ( 345, 184, 138, 22, 1 ) & 4 & ( 15, 8, 6, 1 ) &\cr
& -36 & 102 & ( 170, 85, 68, 16, 1 ) & 4 & ( 10, 5, 4, 1 ) &\cr
\noalign{\hrule}
& -36 & 98 & ( 171, 95, 57, 18, 1 ) & 4 & ( 9, 5, 3, 1 ) &\cr
& -36 & 98 & ( 162, 90, 54, 17, 1 ) & 4 & ( 9, 5, 3, 1 ) &\cr
& -36 & 44 & ( 30, 30, 20, 9, 1 ) & 4 & ( 3, 3, 2, 1 ) &\cr
& -36 & 44 & ( 27, 27, 18, 8, 1 ) & 4 & ( 3, 3, 2, 1 ) &\cr
& -36 & 26 & ( 30, 15, 8, 4, 3 ) & 3 & ( 10, 5, 4, 1 ) &\cr
& -36 & 26 & ( 27, 15, 8, 6, 4 ) & 3 & ( 9, 5, 4, 2 ) &\cr
& -36 & 20 & ( 42, 19, 12, 8, 3 ) & 2 & ( 7, 4, 2, 1 ) &\cr
& -36 & 20 & ( 14, 12, 7, 6, 3 ) & 1 & ( 7, 4, 2, 1 ) &\cr
& -36 & 20 & ( 8, 7, 7, 4, 2 ) & 2 & ( 7, 4, 2, 1 ) &\cr
& -36 & 20 & ( 10, 9, 9, 6, 2 ) & 2 & ( 9, 5, 3, 1 ) &\cr
\noalign{\hrule} }}}}$$
\vfill\eject

$${\fivepoint{
\vbox{\offinterlineskip\tabskip=0pt
\halign{\strut\vrule#
&~~$#$~~\hfil
&~~$#$~~\hfil
&~~$#$~~\hfil
&~~$#$~~\hfil
&~~$#$~~\hfil
&\vrule#\cr
\noalign{\hrule}
&Euler~\#& h^{1,1}& \hat w=(w_1,..,w_5) &i &~~~\bar w&\cr
\noalign{\hrule}
& -36 & 17 & ( 30, 10, 9, 6, 5 ) & 2 & ( 10, 5, 3, 2 ) &\cr
& -36 & 14 & ( 10, 8, 4, 3, 3 ) & 4 & ( 5, 4, 3, 2 ) &\cr
& -40 & 70 & ( 66, 55, 44, 10, 1 ) & 4 & ( 6, 5, 4, 1 ) &\cr
& -40 & 59 & ( 50, 40, 30, 9, 1 ) & 4 & ( 5, 4, 3, 1 ) &\cr
& -40 & 49 & ( 44, 33, 22, 10, 1 ) & 4 & ( 4, 3, 2, 1 ) &\cr
& -40 & 49 & ( 40, 30, 20, 9, 1 ) & 4 & ( 4, 3, 2, 1 ) &\cr
& -40 & 25 & ( 40, 20, 12, 5, 3 ) & 3 & ( 8, 4, 3, 1 ) &\cr
& -40 & 19 & ( 20, 9, 8, 4, 3 ) & 2 & ( 5, 3, 2, 1 ) &\cr
& -42 & 23 & ( 21, 15, 9, 4, 2 ) & 4 & ( 7, 5, 3, 2 ) &\cr
& -42 & 17 & ( 21, 9, 8, 4, 3 ) & 3 & ( 7, 4, 3, 1 ) &\cr
\noalign{\hrule}
& -48 & 77 & ( 84, 72, 48, 11, 1 ) & 4 & ( 7, 6, 4, 1 ) &\cr
& -48 & 67 & ( 165, 110, 33, 20, 2 ) & 4 & ( 15, 10, 3, 2 ) &\cr
& -48 & 67 & ( 66, 55, 33, 10, 1 ) & 4 & ( 6, 5, 3, 1 ) &\cr
& -48 & 59 & ( 65, 52, 26, 12, 1 ) & 4 & ( 5, 4, 2, 1 ) &\cr
& -48 & 59 & ( 60, 48, 24, 11, 1 ) & 4 & ( 5, 4, 2, 1 ) &\cr
& -48 & 50 & ( 36, 36, 27, 8, 1 ) & 4 & ( 4, 4, 3, 1 ) &\cr
& -48 & 43 & ( 126, 84, 31, 7, 4 ) & 3 & ( 9, 6, 2, 1 ) &\cr
& -48 & 39 & ( 24, 24, 16, 7, 1 ) & 4 & ( 3, 3, 2, 1 ) &\cr
& -48 & 35 & ( 40, 28, 9, 4, 3 ) & 3 & ( 10, 7, 3, 1 ) &\cr
& -48 & 35 & ( 27, 18, 18, 8, 1 ) & 4 & ( 3, 2, 2, 1 ) &\cr
\noalign{\hrule}
& -48 & 35 & ( 24, 16, 16, 7, 1 ) & 4 & ( 3, 2, 2, 1 ) &\cr
& -48 & 31 & ( 22, 16, 4, 3, 3 ) & 4 & ( 11, 8, 3, 2 ) &\cr
& -48 & 22 & ( 24, 13, 9, 6, 2 ) & 2 & ( 4, 3, 1, 1 ) &\cr
& -48 & 21 & ( 14, 5, 5, 4, 2 ) & 2 & ( 7, 5, 2, 1 ) &\cr
& -48 & 19 & ( 8, 8, 3, 3, 2 ) & 3 & ( 4, 4, 3, 1 ) &\cr
& -48 & 19 & ( 12, 9, 9, 4, 2 ) & 2 & ( 9, 6, 2, 1 ) &\cr
& -48 & 19 & ( 12, 9, 9, 4, 2 ) & 4 & ( 4, 3, 3, 2 ) &\cr
& -48 & 19 & ( 16, 15, 8, 6, 3 ) & 1 & ( 8, 5, 2, 1 ) &\cr
& -48 & 17 & ( 45, 18, 12, 10, 5 ) & 4 & ( 15, 6, 5, 4 ) &\cr
& -48 & 17 & ( 15, 8, 6, 4, 3 ) & 2 & ( 5, 4, 2, 1 ) &\cr
\noalign{\hrule} }}}}$$
\vfill\eject

$${\fivepoint{
\vbox{\offinterlineskip\tabskip=0pt
\halign{\strut\vrule#
&~~$#$~~\hfil
&~~$#$~~\hfil
&~~$#$~~\hfil
&~~$#$~~\hfil
&~~$#$~~\hfil
&\vrule#\cr
\noalign{\hrule}
&Euler~\#& h^{1,1}& \hat w=(w_1,..,w_5) &i &~~~\bar w&\cr
\noalign{\hrule}
& -48 & 15 & ( 16, 12, 9, 8, 3 ) & 1 & ( 8, 4, 3, 1 ) &\cr
& -48 & 15 & ( 16, 12, 9, 8, 3 ) & 3 & ( 4, 3, 3, 2 ) &\cr
& -48 & 15 & ( 20, 15, 12, 10, 3 ) & 1 & ( 10, 5, 4, 1 ) &\cr
& -48 & 15 & ( 20, 15, 12, 10, 3 ) & 3 & ( 4, 3, 3, 2 ) &\cr
& -48 & 15 & ( 12, 11, 6, 4, 3 ) & 2 & ( 2, 2, 1, 1 ) &\cr
& -48 & 15 & ( 10, 6, 6, 5, 3 ) & 1 & ( 5, 2, 2, 1 ) &\cr
& -48 & 15 & ( 5, 5, 4, 4, 2 ) & 1 & ( 5, 2, 2, 1 ) &\cr
& -48 & 14 & ( 14, 8, 4, 3, 3 ) & 4 & ( 7, 4, 3, 2 ) &\cr
& -48 & 12 & ( 12, 8, 5, 5, 2 ) & 3 & ( 6, 5, 4, 1 ) &\cr
& -48 & 11 & ( 15, 10, 9, 6, 5 ) & 2 & ( 5, 5, 3, 2 ) &\cr
\noalign{\hrule}
& -48 & 11 & ( 20, 15, 12, 8, 5 ) & 2 & ( 5, 5, 3, 2 ) &\cr
& -50 & 24 & ( 25, 15, 8, 5, 2 ) & 3 & ( 5, 3, 2, 1 ) &\cr
& -54 & 53 & ( 45, 36, 27, 8, 1 ) & 4 & ( 5, 4, 3, 1 ) &\cr
& -56 & 93 & ( 160, 80, 64, 15, 1 ) & 4 & ( 10, 5, 4, 1 ) &\cr
& -56 & 76 & ( 112, 56, 42, 13, 1 ) & 4 & ( 8, 4, 3, 1 ) &\cr
& -60 & 222 & ( 777, 518, 222, 36, 1 ) & 4 & ( 21, 14, 6, 1 ) &\cr
& -60 & 164 & ( 465, 310, 124, 30, 1 ) & 4 & ( 15, 10, 4, 1 ) &\cr
& -60 & 164 & ( 450, 300, 120, 29, 1 ) & 4 & ( 15, 10, 4, 1 ) &\cr
& -60 & 19 & ( 25, 20, 15, 12, 3 ) & 4 & ( 5, 4, 3, 3 ) &\cr
& -60 & 14 & ( 9, 9, 4, 3, 2 ) & 3 & ( 3, 3, 2, 1 ) &\cr
\noalign{\hrule}
& -64 & 39 & ( 32, 24, 16, 7, 1 ) & 4 & ( 4, 3, 2, 1 ) &\cr
& -64 & 29 & ( 21, 14, 14, 6, 1 ) & 4 & ( 3, 2, 2, 1 ) &\cr
& -64 & 17 & ( 28, 14, 8, 3, 3 ) & 4 & ( 14, 7, 4, 3 ) &\cr
& -64 & 15 & ( 40, 16, 15, 5, 4 ) & 2 & ( 8, 4, 3, 1 ) &\cr
& -64 & 15 & ( 40, 16, 15, 5, 4 ) & 3 & ( 10, 5, 4, 1 ) &\cr
& -64 & 11 & ( 28, 8, 7, 7, 6 ) & 3 & ( 14, 7, 4, 3 ) &\cr
& -64 & 11 & ( 14, 4, 4, 3, 3 ) & 4 & ( 7, 3, 2, 2 ) &\cr
& -64 & 8 & ( 4, 4, 3, 3, 2 ) & 3 & ( 3, 2, 2, 1 ) &\cr
& -66 & 32 & ( 21, 21, 14, 6, 1 ) & 4 & ( 3, 3, 2, 1 ) &\cr
& -66 & 23 & ( 27, 21, 9, 4, 2 ) & 4 & ( 9, 7, 3, 2 ) &\cr
\noalign{\hrule} }}}}$$
\vfill\eject

$${\fivepoint{
\vbox{\offinterlineskip\tabskip=0pt
\halign{\strut\vrule#
&~~$#$~~\hfil
&~~$#$~~\hfil
&~~$#$~~\hfil
&~~$#$~~\hfil
&~~$#$~~\hfil
&\vrule#\cr
\noalign{\hrule}
&Euler~\#& h^{1,1}& \hat w=(w_1,..,w_5) &i &~~~\bar w&\cr
\noalign{\hrule}
& -72 & 88 & ( 117, 104, 78, 12, 1 ) & 4 & ( 9, 8, 6, 1 ) &\cr
& -72 & 69 & ( 104, 52, 39, 12, 1 ) & 4 & ( 8, 4, 3, 1 ) &\cr
& -72 & 57 & ( 60, 50, 30, 9, 1 ) & 4 & ( 6, 5, 3, 1 ) &\cr
& -72 & 56 & ( 54, 45, 36, 8, 1 ) & 4 & ( 6, 5, 4, 1 ) &\cr
& -72 & 53 & ( 78, 39, 26, 12, 1 ) & 4 & ( 6, 3, 2, 1 ) &\cr
& -72 & 53 & ( 72, 36, 24, 11, 1 ) & 4 & ( 6, 3, 2, 1 ) &\cr
& -72 & 40 & ( 28, 28, 21, 6, 1 ) & 4 & ( 4, 4, 3, 1 ) &\cr
& -72 & 32 & ( 60, 35, 15, 8, 2 ) & 4 & ( 12, 7, 3, 2 ) &\cr
& -72 & 29 & ( 42, 27, 8, 4, 3 ) & 3 & ( 14, 9, 4, 1 ) &\cr
& -72 & 26 & ( 40, 20, 13, 5, 2 ) & 3 & ( 4, 2, 1, 1 ) &\cr
\noalign{\hrule}
& -72 & 26 & ( 18, 12, 12, 5, 1 ) & 4 & ( 3, 2, 2, 1 ) &\cr
& -72 & 23 & ( 42, 21, 12, 7, 2 ) & 3 & ( 6, 3, 2, 1 ) &\cr
& -72 & 23 & ( 35, 21, 14, 12, 2 ) & 4 & ( 5, 3, 2, 2 ) &\cr
& -72 & 21 & ( 14, 8, 3, 3, 2 ) & 3 & ( 7, 4, 3, 1 ) &\cr
& -72 & 20 & ( 36, 17, 12, 4, 3 ) & 2 & ( 3, 1, 1, 1 ) &\cr
& -72 & 14 & ( 20, 10, 4, 3, 3 ) & 4 & ( 10, 5, 3, 2 ) &\cr
& -72 & 14 & ( 16, 9, 4, 4, 3 ) & 2 & ( 4, 3, 1, 1 ) &\cr
& -72 & 13 & ( 30, 12, 10, 5, 3 ) & 2 & ( 6, 3, 2, 1 ) &\cr
& -72 & 13 & ( 30, 12, 10, 5, 3 ) & 3 & ( 10, 5, 4, 1 ) &\cr
& -72 & 13 & ( 24, 9, 8, 4, 3 ) & 2 & ( 6, 3, 2, 1 ) &\cr
\noalign{\hrule}
& -72 & 13 & ( 24, 9, 8, 4, 3 ) & 3 & ( 8, 4, 3, 1 ) &\cr
& -72 & 13 & ( 24, 11, 6, 4, 3 ) & 2 & ( 4, 2, 1, 1 ) &\cr
& -72 & 13 & ( 30, 12, 8, 5, 5 ) & 4 & ( 15, 6, 5, 4 ) &\cr
& -72 & 13 & ( 12, 7, 6, 3, 2 ) & 2 & ( 2, 1, 1, 1 ) &\cr
& -72 & 11 & ( 7, 6, 6, 3, 2 ) & 1 & ( 1, 1, 1, 1 ) &\cr
& -72 & 10 & ( 20, 6, 5, 5, 4 ) & 3 & ( 10, 5, 3, 2 ) &\cr
& -72 & 9 & ( 8, 4, 3, 3, 2 ) & 3 & ( 4, 3, 2, 1 ) &\cr
& -72 & 8 & ( 12, 9, 8, 4, 3 ) & 2 & ( 3, 3, 2, 1 ) &\cr
& -72 & 8 & ( 12, 9, 8, 4, 3 ) & 3 & ( 4, 4, 3, 1 ) &\cr
& -72 & 7 & ( 10, 6, 5, 5, 4 ) & 3 & ( 5, 5, 3, 2 ) &\cr
\noalign{\hrule} }}}}$$
\vfill\eject

$${\fivepoint{
\vbox{\offinterlineskip\tabskip=0pt
\halign{\strut\vrule#
&~~$#$~~\hfil
&~~$#$~~\hfil
&~~$#$~~\hfil
&~~$#$~~\hfil
&~~$#$~~\hfil
&\vrule#\cr
\noalign{\hrule}
&Euler~\#& h^{1,1}& \hat w=(w_1,..,w_5) &i &~~~\bar w&\cr
\noalign{\hrule}
& -80 & 39 & ( 55, 22, 22, 10, 1 ) & 4 & ( 5, 2, 2, 1 ) &\cr
& -80 & 39 & ( 50, 20, 20, 9, 1 ) & 4 & ( 5, 2, 2, 1 ) &\cr
& -80 & 33 & ( 28, 21, 14, 6, 1 ) & 4 & ( 4, 3, 2, 1 ) &\cr
& -80 & 11 & ( 15, 10, 8, 5, 2 ) & 3 & ( 3, 2, 2, 1 ) &\cr
& -84 & 76 & ( 126, 70, 42, 13, 1 ) & 4 & ( 9, 5, 3, 1 ) &\cr
& -84 & 62 & ( 105, 60, 30, 14, 1 ) & 4 & ( 7, 4, 2, 1 ) &\cr
& -84 & 62 & ( 98, 56, 28, 13, 1 ) & 4 & ( 7, 4, 2, 1 ) &\cr
& -84 & 61 & ( 70, 60, 40, 9, 1 ) & 4 & ( 7, 6, 4, 1 ) &\cr
& -84 & 48 & ( 66, 33, 22, 10, 1 ) & 4 & ( 6, 3, 2, 1 ) &\cr
& -84 & 44 & ( 45, 36, 18, 8, 1 ) & 4 & ( 5, 4, 2, 1 ) &\cr
\noalign{\hrule}
& -84 & 41 & ( 35, 28, 21, 6, 1 ) & 4 & ( 5, 4, 3, 1 ) &\cr
& -84 & 14 & ( 15, 9, 4, 3, 2 ) & 3 & ( 5, 3, 2, 1 ) &\cr
& -84 & 12 & ( 10, 9, 5, 3, 3 ) & 1 & ( 5, 3, 1, 1 ) &\cr
& -84 & 12 & ( 6, 5, 5, 2, 2 ) & 2 & ( 5, 3, 1, 1 ) &\cr
& -84 & 9 & ( 10, 4, 3, 3, 2 ) & 3 & ( 5, 3, 2, 1 ) &\cr
& -84 & 8 & ( 9, 8, 4, 3, 3 ) & 2 & ( 4, 3, 1, 1 ) &\cr
& -88 & 17 & ( 34, 20, 8, 3, 3 ) & 4 & ( 17, 10, 4, 3 ) &\cr
& -90 & 24 & ( 15, 15, 10, 4, 1 ) & 4 & ( 3, 3, 2, 1 ) &\cr
& -90 & 24 & ( 35, 25, 8, 5, 2 ) & 3 & ( 7, 5, 2, 1 ) &\cr
& -90 & 8 & ( 9, 4, 3, 3, 2 ) & 2 & ( 3, 2, 1, 1 ) &\cr
\noalign{\hrule}
& -92 & 5 & ( 4, 3, 3, 2, 2 ) & 2 & ( 3, 2, 1, 1 ) &\cr
& -96 & 147 & ( 405, 270, 108, 26, 1 ) & 4 & ( 15, 10, 4, 1 ) &\cr
& -96 & 119 & ( 300, 200, 75, 24, 1 ) & 4 & ( 12, 8, 3, 1 ) &\cr
& -96 & 119 & ( 288, 192, 72, 23, 1 ) & 4 & ( 12, 8, 3, 1 ) &\cr
& -96 & 47 & ( 48, 40, 24, 7, 1 ) & 4 & ( 6, 5, 3, 1 ) &\cr
& -96 & 43 & ( 60, 30, 20, 9, 1 ) & 4 & ( 6, 3, 2, 1 ) &\cr
& -96 & 39 & ( 40, 32, 16, 7, 1 ) & 4 & ( 5, 4, 2, 1 ) &\cr
& -96 & 33 & ( 45, 18, 18, 8, 1 ) & 4 & ( 5, 2, 2, 1 ) &\cr
& -96 & 31 & ( 96, 64, 21, 8, 3 ) & 3 & ( 12, 8, 3, 1 ) &\cr
& -96 & 31 & ( 72, 45, 16, 9, 2 ) & 3 & ( 8, 5, 2, 1 ) &\cr
\noalign{\hrule} }}}}$$
\vfill\eject

$${\fivepoint{
\vbox{\offinterlineskip\tabskip=0pt
\halign{\strut\vrule#
&~~$#$~~\hfil
&~~$#$~~\hfil
&~~$#$~~\hfil
&~~$#$~~\hfil
&~~$#$~~\hfil
&\vrule#\cr
\noalign{\hrule}
&Euler~\#& h^{1,1}& \hat w=(w_1,..,w_5) &i &~~~\bar w&\cr
\noalign{\hrule}
& -96 & 27 & ( 24, 18, 12, 5, 1 ) & 4 & ( 4, 3, 2, 1 ) &\cr
& -96 & 23 & ( 60, 40, 9, 8, 3 ) & 3 & ( 15, 10, 3, 2 ) &\cr
& -96 & 19 & ( 30, 15, 9, 4, 2 ) & 4 & ( 10, 5, 3, 2 ) &\cr
& -96 & 19 & ( 15, 10, 10, 4, 1 ) & 4 & ( 3, 2, 2, 1 ) &\cr
& -96 & 17 & ( 12, 12, 6, 5, 1 ) & 4 & ( 2, 2, 1, 1 ) &\cr
& -96 & 17 & ( 14, 14, 7, 6, 1 ) & 4 & ( 2, 2, 1, 1 ) &\cr
& -96 & 14 & ( 24, 10, 8, 3, 3 ) & 4 & ( 12, 5, 4, 3 ) &\cr
& -96 & 11 & ( 18, 8, 4, 3, 3 ) & 2 & ( 6, 4, 1, 1 ) &\cr
& -96 & 11 & ( 18, 8, 4, 3, 3 ) & 4 & ( 9, 4, 3, 2 ) &\cr
& -96 & 11 & ( 8, 3, 3, 2, 2 ) & 2 & ( 4, 3, 1, 1 ) &\cr
\noalign{\hrule}
& -96 & 9 & ( 20, 9, 4, 4, 3 ) & 2 & ( 5, 3, 1, 1 ) &\cr
& -96 & 7 & ( 8, 6, 4, 3, 3 ) & 1 & ( 4, 2, 1, 1 ) &\cr
& -96 & 7 & ( 8, 6, 4, 3, 3 ) & 4 & ( 4, 3, 3, 2 ) &\cr
& -96 & 7 & ( 9, 6, 4, 3, 2 ) & 3 & ( 3, 2, 2, 1 ) &\cr
& -96 & 5 & ( 6, 4, 3, 3, 2 ) & 2 & ( 2, 2, 1, 1 ) &\cr
& -96 & 5 & ( 6, 4, 3, 3, 2 ) & 3 & ( 3, 3, 2, 1 ) &\cr
& -96 & 5 & ( 4, 3, 3, 3, 2 ) & 1 & ( 2, 1, 1, 1 ) &\cr
& -104 & 17 & ( 12, 8, 8, 3, 1 ) & 4 & ( 3, 2, 2, 1 ) &\cr
& -108 & 14 & ( 26, 16, 4, 3, 3 ) & 4 & ( 13, 8, 3, 2 ) &\cr
& -108 & 13 & ( 10, 10, 5, 4, 1 ) & 4 & ( 2, 2, 1, 1 ) &\cr
\noalign{\hrule}
& -108 & 6 & ( 3, 3, 2, 2, 2 ) & 1 & ( 3, 1, 1, 1 ) &\cr
& -112 & 20 & ( 49, 21, 14, 12, 2 ) & 4 & ( 7, 3, 2, 2 ) &\cr
& -112 & 20 & ( 21, 14, 7, 6, 1 ) & 4 & ( 3, 2, 1, 1 ) &\cr
& -112 & 20 & ( 24, 16, 8, 7, 1 ) & 4 & ( 3, 2, 1, 1 ) &\cr
& -112 & 10 & ( 16, 8, 3, 3, 2 ) & 3 & ( 8, 4, 3, 1 ) &\cr
& -112 & 7 & ( 20, 8, 5, 5, 2 ) & 2 & ( 4, 2, 1, 1 ) &\cr
& -112 & 7 & ( 20, 8, 5, 5, 2 ) & 3 & ( 10, 5, 4, 1 ) &\cr
& -120 & 109 & ( 240, 128, 96, 15, 1 ) & 4 & ( 15, 8, 6, 1 ) &\cr
& -120 & 108 & ( 264, 176, 66, 21, 1 ) & 4 & ( 12, 8, 3, 1 ) &\cr
& -120 & 76 & ( 132, 72, 48, 11, 1 ) & 4 & ( 11, 6, 4, 1 ) &\cr
\noalign{\hrule} }}}}$$
\vfill\eject

$${\fivepoint{
\vbox{\offinterlineskip\tabskip=0pt
\halign{\strut\vrule#
&~~$#$~~\hfil
&~~$#$~~\hfil
&~~$#$~~\hfil
&~~$#$~~\hfil
&~~$#$~~\hfil
&\vrule#\cr
\noalign{\hrule}
&Euler~\#& h^{1,1}& \hat w=(w_1,..,w_5) &i &~~~\bar w&\cr
\noalign{\hrule}
& -120 & 68 & ( 81, 72, 54, 8, 1 ) & 4 & ( 9, 8, 6, 1 ) &\cr
& -120 & 65 & ( 110, 55, 44, 10, 1 ) & 4 & ( 10, 5, 4, 1 ) &\cr
& -120 & 49 & ( 72, 36, 27, 8, 1 ) & 4 & ( 8, 4, 3, 1 ) &\cr
& -120 & 48 & ( 49, 42, 28, 6, 1 ) & 4 & ( 7, 6, 4, 1 ) &\cr
& -120 & 47 & ( 77, 44, 22, 10, 1 ) & 4 & ( 7, 4, 2, 1 ) &\cr
& -120 & 39 & ( 36, 30, 18, 5, 1 ) & 4 & ( 6, 5, 3, 1 ) &\cr
& -120 & 38 & ( 36, 30, 24, 5, 1 ) & 4 & ( 6, 5, 4, 1 ) &\cr
& -120 & 31 & ( 35, 28, 14, 6, 1 ) & 4 & ( 5, 4, 2, 1 ) &\cr
& -120 & 29 & ( 25, 20, 15, 4, 1 ) & 4 & ( 5, 4, 3, 1 ) &\cr
& -120 & 26 & ( 20, 20, 15, 4, 1 ) & 4 & ( 4, 4, 3, 1 ) &\cr
\noalign{\hrule}
& -120 & 26 & ( 16, 16, 12, 3, 1 ) & 4 & ( 4, 4, 3, 1 ) &\cr
& -120 & 25 & ( 30, 12, 12, 5, 1 ) & 4 & ( 5, 2, 2, 1 ) &\cr
& -120 & 25 & ( 20, 14, 3, 3, 2 ) & 3 & ( 10, 7, 3, 1 ) &\cr
& -120 & 22 & ( 60, 40, 12, 5, 3 ) & 3 & ( 12, 8, 3, 1 ) &\cr
& -120 & 22 & ( 40, 25, 8, 5, 2 ) & 3 & ( 8, 5, 2, 1 ) &\cr
& -120 & 21 & ( 36, 21, 9, 4, 2 ) & 4 & ( 12, 7, 3, 2 ) &\cr
& -120 & 17 & ( 12, 12, 8, 3, 1 ) & 4 & ( 3, 3, 2, 1 ) &\cr
& -120 & 17 & ( 18, 12, 6, 5, 1 ) & 4 & ( 3, 2, 1, 1 ) &\cr
& -120 & 10 & ( 18, 9, 4, 3, 2 ) & 3 & ( 6, 3, 2, 1 ) &\cr
& -120 & 10 & ( 8, 8, 4, 3, 1 ) & 4 & ( 2, 2, 1, 1 ) &\cr
\noalign{\hrule}
& -120 & 10 & ( 15, 9, 6, 4, 2 ) & 4 & ( 5, 3, 2, 2 ) &\cr
& -120 & 9 & ( 25, 10, 8, 5, 2 ) & 3 & ( 5, 2, 2, 1 ) &\cr
& -120 & 9 & ( 10, 5, 5, 4, 1 ) & 4 & ( 2, 1, 1, 1 ) &\cr
& -120 & 9 & ( 12, 6, 6, 5, 1 ) & 4 & ( 2, 1, 1, 1 ) &\cr
& -120 & 6 & ( 12, 4, 3, 3, 2 ) & 2 & ( 4, 2, 1, 1 ) &\cr
& -120 & 6 & ( 12, 4, 3, 3, 2 ) & 3 & ( 6, 3, 2, 1 ) &\cr
& -120 & 5 & ( 4, 4, 4, 3, 1 ) & 4 & ( 1, 1, 1, 1 ) &\cr
& -120 & 5 & ( 5, 5, 5, 4, 1 ) & 4 & ( 1, 1, 1, 1 ) &\cr
& -128 & 19 & ( 16, 12, 8, 3, 1 ) & 4 & ( 4, 3, 2, 1 ) &\cr
& -128 & 7 & ( 8, 4, 4, 3, 1 ) & 4 & ( 2, 1, 1, 1 ) &\cr
\noalign{\hrule} }}}}$$
\vfill\eject

$${\fivepoint{
\vbox{\offinterlineskip\tabskip=0pt
\halign{\strut\vrule#
&~~$#$~~\hfil
&~~$#$~~\hfil
&~~$#$~~\hfil
&~~$#$~~\hfil
&~~$#$~~\hfil
&\vrule#\cr
\noalign{\hrule}
&Euler~\#& h^{1,1}& \hat w=(w_1,..,w_5) &i &~~~\bar w&\cr
\noalign{\hrule}
& -130 & 14 & ( 15, 10, 5, 4, 1 ) & 4 & ( 3, 2, 1, 1 ) &\cr
& -132 & 56 & ( 90, 50, 30, 9, 1 ) & 4 & ( 9, 5, 3, 1 ) &\cr
& -132 & 30 & ( 42, 21, 14, 6, 1 ) & 4 & ( 6, 3, 2, 1 ) &\cr
& -132 & 14 & ( 21, 15, 4, 3, 2 ) & 3 & ( 7, 5, 2, 1 ) &\cr
& -132 & 7 & ( 10, 8, 6, 3, 3 ) & 4 & ( 5, 4, 3, 3 ) &\cr
& -132 & 4 & ( 10, 3, 3, 2, 2 ) & 2 & ( 5, 3, 1, 1 ) &\cr
& -132 & 3 & ( 3, 3, 3, 2, 1 ) & 4 & ( 1, 1, 1, 1 ) &\cr
& -136 & 34 & ( 30, 25, 20, 4, 1 ) & 4 & ( 6, 5, 4, 1 ) &\cr
& -138 & 14 & ( 9, 9, 6, 2, 1 ) & 4 & ( 3, 3, 2, 1 ) &\cr
& -144 & 71 & ( 171, 114, 38, 18, 1 ) & 4 & ( 9, 6, 2, 1 ) &\cr
\noalign{\hrule}
& -144 & 71 & ( 162, 108, 36, 17, 1 ) & 4 & ( 9, 6, 2, 1 ) &\cr
& -144 & 55 & ( 90, 45, 36, 8, 1 ) & 4 & ( 10, 5, 4, 1 ) &\cr
& -144 & 38 & ( 56, 32, 16, 7, 1 ) & 4 & ( 7, 4, 2, 1 ) &\cr
& -144 & 26 & ( 36, 18, 12, 5, 1 ) & 4 & ( 6, 3, 2, 1 ) &\cr
& -144 & 26 & ( 81, 54, 16, 9, 2 ) & 3 & ( 9, 6, 2, 1 ) &\cr
& -144 & 26 & ( 36, 27, 9, 8, 1 ) & 4 & ( 4, 3, 1, 1 ) &\cr
& -144 & 26 & ( 40, 30, 10, 9, 1 ) & 4 & ( 4, 3, 1, 1 ) &\cr
& -144 & 25 & ( 25, 20, 10, 4, 1 ) & 4 & ( 5, 4, 2, 1 ) &\cr
& -144 & 23 & ( 20, 16, 12, 3, 1 ) & 4 & ( 5, 4, 3, 1 ) &\cr
& -144 & 19 & ( 32, 16, 8, 7, 1 ) & 4 & ( 4, 2, 1, 1 ) &\cr
\noalign{\hrule}
& -144 & 19 & ( 36, 18, 9, 8, 1 ) & 4 & ( 4, 2, 1, 1 ) &\cr
& -144 & 11 & ( 12, 8, 4, 3, 1 ) & 4 & ( 3, 2, 1, 1 ) &\cr
& -144 & 11 & ( 9, 6, 6, 2, 1 ) & 4 & ( 3, 2, 2, 1 ) &\cr
& -144 & 7 & ( 6, 6, 3, 2, 1 ) & 4 & ( 2, 2, 1, 1 ) &\cr
& -144 & 5 & ( 15, 6, 4, 3, 2 ) & 3 & ( 5, 2, 2, 1 ) &\cr
& -144 & 5 & ( 6, 3, 3, 2, 1 ) & 4 & ( 2, 1, 1, 1 ) &\cr
& -150 & 29 & ( 30, 25, 15, 4, 1 ) & 4 & ( 6, 5, 3, 1 ) &\cr
& -152 & 37 & ( 56, 28, 21, 6, 1 ) & 4 & ( 8, 4, 3, 1 ) &\cr
& -152 & 17 & ( 25, 10, 10, 4, 1 ) & 4 & ( 5, 2, 2, 1 ) &\cr
& -152 & 16 & ( 28, 14, 7, 6, 1 ) & 4 & ( 4, 2, 1, 1 ) &\cr
\noalign{\hrule} }}}}$$
\vfill\eject

$${\fivepoint{
\vbox{\offinterlineskip\tabskip=0pt
\halign{\strut\vrule#
&~~$#$~~\hfil
&~~$#$~~\hfil
&~~$#$~~\hfil
&~~$#$~~\hfil
&~~$#$~~\hfil
&\vrule#\cr
\noalign{\hrule}
&Euler~\#& h^{1,1}& \hat w=(w_1,..,w_5) &i &~~~\bar w&\cr
\noalign{\hrule}
& -156 & 176 & ( 609, 406, 174, 28, 1 ) & 4 & ( 21, 14, 6, 1 ) &\cr
& -156 & 66 & ( 153, 102, 34, 16, 1 ) & 4 & ( 9, 6, 2, 1 ) &\cr
& -156 & 21 & ( 63, 42, 12, 7, 2 ) & 3 & ( 9, 6, 2, 1 ) &\cr
& -156 & 21 & ( 28, 21, 7, 6, 1 ) & 4 & ( 4, 3, 1, 1 ) &\cr
& -156 & 8 & ( 18, 6, 6, 5, 1 ) & 4 & ( 3, 1, 1, 1 ) &\cr
& -156 & 8 & ( 21, 7, 7, 6, 1 ) & 4 & ( 3, 1, 1, 1 ) &\cr
& -160 & 59 & ( 110, 60, 40, 9, 1 ) & 4 & ( 11, 6, 4, 1 ) &\cr
& -160 & 15 & ( 20, 8, 8, 3, 1 ) & 4 & ( 5, 2, 2, 1 ) &\cr
& -168 & 86 & ( 204, 136, 51, 16, 1 ) & 4 & ( 12, 8, 3, 1 ) &\cr
& -168 & 50 & ( 63, 56, 42, 6, 1 ) & 4 & ( 9, 8, 6, 1 ) &\cr
\noalign{\hrule}
& -168 & 32 & ( 48, 24, 18, 5, 1 ) & 4 & ( 8, 4, 3, 1 ) &\cr
& -168 & 25 & ( 24, 20, 16, 3, 1 ) & 4 & ( 6, 5, 4, 1 ) &\cr
& -168 & 20 & ( 30, 15, 10, 4, 1 ) & 4 & ( 6, 3, 2, 1 ) &\cr
& -168 & 18 & ( 20, 16, 8, 3, 1 ) & 4 & ( 5, 4, 2, 1 ) &\cr
& -168 & 17 & ( 24, 18, 6, 5, 1 ) & 4 & ( 4, 3, 1, 1 ) &\cr
& -168 & 16 & ( 12, 12, 9, 2, 1 ) & 4 & ( 4, 4, 3, 1 ) &\cr
& -168 & 13 & ( 12, 9, 6, 2, 1 ) & 4 & ( 4, 3, 2, 1 ) &\cr
& -168 & 12 & ( 20, 10, 5, 4, 1 ) & 4 & ( 4, 2, 1, 1 ) &\cr
& -168 & 12 & ( 24, 15, 4, 3, 2 ) & 3 & ( 8, 5, 2, 1 ) &\cr
& -168 & 8 & ( 21, 9, 6, 4, 2 ) & 4 & ( 7, 3, 2, 2 ) &\cr
\noalign{\hrule}
& -168 & 8 & ( 9, 6, 3, 2, 1 ) & 4 & ( 3, 2, 1, 1 ) &\cr
& -168 & 6 & ( 12, 4, 4, 3, 1 ) & 4 & ( 3, 1, 1, 1 ) &\cr
& -168 & 2 & ( 2, 2, 2, 1, 1 ) & 4 & ( 1, 1, 1, 1 ) &\cr
& -172 & 29 & ( 49, 28, 14, 6, 1 ) & 4 & ( 7, 4, 2, 1 ) &\cr
& -172 & 10 & ( 22, 14, 3, 3, 2 ) & 3 & ( 11, 7, 3, 1 ) &\cr
& -180 & 26 & ( 42, 24, 12, 5, 1 ) & 4 & ( 7, 4, 2, 1 ) &\cr
& -180 & 24 & ( 50, 30, 10, 9, 1 ) & 4 & ( 5, 3, 1, 1 ) &\cr
& -180 & 24 & ( 55, 33, 11, 10, 1 ) & 4 & ( 5, 3, 1, 1 ) &\cr
& -180 & 17 & ( 15, 12, 9, 2, 1 ) & 4 & ( 5, 4, 3, 1 ) &\cr
& -180 & 14 & ( 45, 30, 8, 5, 2 ) & 3 & ( 9, 6, 2, 1 ) &\cr
\noalign{\hrule} }}}}$$
\vfill\eject

$${\fivepoint{
\vbox{\offinterlineskip\tabskip=0pt
\halign{\strut\vrule#
&~~$#$~~\hfil
&~~$#$~~\hfil
&~~$#$~~\hfil
&~~$#$~~\hfil
&~~$#$~~\hfil
&\vrule#\cr
\noalign{\hrule}
&Euler~\#& h^{1,1}& \hat w=(w_1,..,w_5) &i &~~~\bar w&\cr
\noalign{\hrule}
& -180 & 14 & ( 20, 15, 5, 4, 1 ) & 4 & ( 4, 3, 1, 1 ) &\cr
& -184 & 9 & ( 16, 8, 4, 3, 1 ) & 4 & ( 4, 2, 1, 1 ) &\cr
& -192 & 77 & ( 180, 96, 72, 11, 1 ) & 4 & ( 15, 8, 6, 1 ) &\cr
& -192 & 51 & ( 117, 78, 26, 12, 1 ) & 4 & ( 9, 6, 2, 1 ) &\cr
& -192 & 47 & ( 88, 48, 32, 7, 1 ) & 4 & ( 11, 6, 4, 1 ) &\cr
& -192 & 41 & ( 54, 48, 36, 5, 1 ) & 4 & ( 9, 8, 6, 1 ) &\cr
& -192 & 35 & ( 63, 35, 21, 6, 1 ) & 4 & ( 9, 5, 3, 1 ) &\cr
& -192 & 27 & ( 28, 24, 16, 3, 1 ) & 4 & ( 7, 6, 4, 1 ) &\cr
& -192 & 21 & ( 24, 20, 12, 3, 1 ) & 4 & ( 6, 5, 3, 1 ) &\cr
& -192 & 19 & ( 40, 24, 8, 7, 1 ) & 4 & ( 5, 3, 1, 1 ) &\cr
\noalign{\hrule}
& -192 & 16 & ( 24, 12, 8, 3, 1 ) & 4 & ( 6, 3, 2, 1 ) &\cr
& -192 & 13 & ( 16, 12, 4, 3, 1 ) & 4 & ( 4, 3, 1, 1 ) &\cr
& -192 & 11 & ( 30, 20, 4, 3, 3 ) & 4 & ( 15, 10, 3, 2 ) &\cr
& -192 & 11 & ( 36, 24, 7, 3, 2 ) & 3 & ( 6, 4, 1, 1 ) &\cr
& -192 & 8 & ( 6, 4, 4, 1, 1 ) & 4 & ( 3, 2, 2, 1 ) &\cr
& -192 & 5 & ( 4, 4, 2, 1, 1 ) & 4 & ( 2, 2, 1, 1 ) &\cr
& -192 & 3 & ( 9, 3, 3, 2, 1 ) & 4 & ( 3, 1, 1, 1 ) &\cr
& -192 & 3 & ( 4, 2, 2, 1, 1 ) & 4 & ( 2, 1, 1, 1 ) &\cr
& -200 & 25 & ( 40, 20, 15, 4, 1 ) & 4 & ( 8, 4, 3, 1 ) &\cr
& -204 & 16 & ( 30, 18, 6, 5, 1 ) & 4 & ( 5, 3, 1, 1 ) &\cr
\noalign{\hrule}
& -204 & 14 & ( 42, 28, 8, 3, 3 ) & 4 & ( 21, 14, 4, 3 ) &\cr
& -204 & 9 & ( 6, 6, 4, 1, 1 ) & 4 & ( 3, 3, 2, 1 ) &\cr
& -216 & 92 & ( 240, 160, 64, 15, 1 ) & 4 & ( 15, 10, 4, 1 ) &\cr
& -216 & 68 & ( 165, 88, 66, 10, 1 ) & 4 & ( 15, 8, 6, 1 ) &\cr
& -216 & 66 & ( 156, 104, 39, 12, 1 ) & 4 & ( 12, 8, 3, 1 ) &\cr
& -216 & 42 & ( 90, 60, 20, 9, 1 ) & 4 & ( 9, 6, 2, 1 ) &\cr
& -216 & 36 & ( 45, 40, 30, 4, 1 ) & 4 & ( 9, 8, 6, 1 ) &\cr
& -216 & 33 & ( 60, 30, 24, 5, 1 ) & 4 & ( 10, 5, 4, 1 ) &\cr
& -216 & 33 & ( 50, 25, 20, 4, 1 ) & 4 & ( 10, 5, 4, 1 ) &\cr
& -216 & 21 & ( 32, 16, 12, 3, 1 ) & 4 & ( 8, 4, 3, 1 ) &\cr
\noalign{\hrule} }}}}$$
\vfill\eject

$${\fivepoint{
\vbox{\offinterlineskip\tabskip=0pt
\halign{\strut\vrule#
&~~$#$~~\hfil
&~~$#$~~\hfil
&~~$#$~~\hfil
&~~$#$~~\hfil
&~~$#$~~\hfil
&\vrule#\cr
\noalign{\hrule}
&Euler~\#& h^{1,1}& \hat w=(w_1,..,w_5) &i &~~~\bar w&\cr
\noalign{\hrule}
& -216 & 18 & ( 18, 15, 12, 2, 1 ) & 4 & ( 6, 5, 4, 1 ) &\cr
& -216 & 17 & ( 45, 30, 9, 4, 2 ) & 4 & ( 15, 10, 3, 2 ) &\cr
& -216 & 17 & ( 18, 15, 9, 2, 1 ) & 4 & ( 6, 5, 3, 1 ) &\cr
& -216 & 13 & ( 15, 12, 6, 2, 1 ) & 4 & ( 5, 4, 2, 1 ) &\cr
& -216 & 9 & ( 15, 6, 6, 2, 1 ) & 4 & ( 5, 2, 2, 1 ) &\cr
& -216 & 6 & ( 12, 6, 3, 2, 1 ) & 4 & ( 4, 2, 1, 1 ) &\cr
& -224 & 17 & ( 28, 16, 8, 3, 1 ) & 4 & ( 7, 4, 2, 1 ) &\cr
& -228 & 12 & ( 18, 9, 6, 2, 1 ) & 4 & ( 6, 3, 2, 1 ) &\cr
& -228 & 8 & ( 27, 18, 4, 3, 2 ) & 3 & ( 9, 6, 2, 1 ) &\cr
& -228 & 8 & ( 12, 9, 3, 2, 1 ) & 4 & ( 4, 3, 1, 1 ) &\cr
\noalign{\hrule}
& -232 & 9 & ( 8, 6, 4, 1, 1 ) & 4 & ( 4, 3, 2, 1 ) &\cr
& -232 & 5 & ( 6, 4, 2, 1, 1 ) & 4 & ( 3, 2, 1, 1 ) &\cr
& -240 & 137 & ( 462, 308, 132, 21, 1 ) & 4 & ( 21, 14, 6, 1 ) &\cr
& -240 & 34 & ( 81, 54, 18, 8, 1 ) & 4 & ( 9, 6, 2, 1 ) &\cr
& -240 & 23 & ( 72, 48, 12, 11, 1 ) & 4 & ( 6, 4, 1, 1 ) &\cr
& -240 & 23 & ( 78, 52, 13, 12, 1 ) & 4 & ( 6, 4, 1, 1 ) &\cr
& -240 & 11 & ( 8, 8, 6, 1, 1 ) & 4 & ( 4, 4, 3, 1 ) &\cr
& -240 & 9 & ( 20, 12, 4, 3, 1 ) & 4 & ( 5, 3, 1, 1 ) &\cr
& -240 & 7 & ( 24, 16, 3, 3, 2 ) & 3 & ( 12, 8, 3, 1 ) &\cr
& -252 & 76 & ( 210, 140, 56, 13, 1 ) & 4 & ( 15, 10, 4, 1 ) &\cr
\noalign{\hrule}
& -252 & 18 & ( 54, 36, 9, 8, 1 ) & 4 & ( 6, 4, 1, 1 ) &\cr
& -252 & 18 & ( 21, 18, 12, 2, 1 ) & 4 & ( 7, 6, 4, 1 ) &\cr
& -252 & 2 & ( 6, 2, 2, 1, 1 ) & 4 & ( 3, 1, 1, 1 ) &\cr
& -256 & 23 & ( 40, 20, 16, 3, 1 ) & 4 & ( 10, 5, 4, 1 ) &\cr
& -264 & 48 & ( 108, 72, 27, 8, 1 ) & 4 & ( 12, 8, 3, 1 ) &\cr
& -264 & 28 & ( 63, 42, 14, 6, 1 ) & 4 & ( 9, 6, 2, 1 ) &\cr
& -264 & 26 & ( 36, 32, 24, 3, 1 ) & 4 & ( 9, 8, 6, 1 ) &\cr
& -264 & 20 & ( 36, 20, 12, 3, 1 ) & 4 & ( 9, 5, 3, 1 ) &\cr
& -264 & 15 & ( 24, 12, 9, 2, 1 ) & 4 & ( 8, 4, 3, 1 ) &\cr
& -264 & 15 & ( 42, 28, 7, 6, 1 ) & 4 & ( 6, 4, 1, 1 ) &\cr
\noalign{\hrule} }}}}$$
\vfill\eject

$${\fivepoint{
\vbox{\offinterlineskip\tabskip=0pt
\halign{\strut\vrule#
&~~$#$~~\hfil
&~~$#$~~\hfil
&~~$#$~~\hfil
&~~$#$~~\hfil
&~~$#$~~\hfil
&\vrule#\cr
\noalign{\hrule}
&Euler~\#& h^{1,1}& \hat w=(w_1,..,w_5) &i &~~~\bar w&\cr
\noalign{\hrule}
& -264 & 11 & ( 10, 8, 6, 1, 1 ) & 4 & ( 5, 4, 3, 1 ) &\cr
& -272 & 7 & ( 10, 4, 4, 1, 1 ) & 4 & ( 5, 2, 2, 1 ) &\cr
& -276 & 48 & ( 105, 56, 42, 6, 1 ) & 4 & ( 15, 8, 6, 1 ) &\cr
& -276 & 6 & ( 15, 9, 3, 2, 1 ) & 4 & ( 5, 3, 1, 1 ) &\cr
& -288 & 115 & ( 399, 266, 114, 18, 1 ) & 4 & ( 21, 14, 6, 1 ) &\cr
& -288 & 63 & ( 165, 110, 44, 10, 1 ) & 4 & ( 15, 10, 4, 1 ) &\cr
& -288 & 39 & ( 96, 64, 24, 7, 1 ) & 4 & ( 12, 8, 3, 1 ) &\cr
& -288 & 23 & ( 44, 24, 16, 3, 1 ) & 4 & ( 11, 6, 4, 1 ) &\cr
& -288 & 11 & ( 30, 20, 5, 4, 1 ) & 4 & ( 6, 4, 1, 1 ) &\cr
& -288 & 11 & ( 21, 12, 6, 2, 1 ) & 4 & ( 7, 4, 2, 1 ) &\cr
\noalign{\hrule}
& -288 & 9 & ( 10, 8, 4, 1, 1 ) & 4 & ( 5, 4, 2, 1 ) &\cr
& -288 & 4 & ( 8, 4, 2, 1, 1 ) & 4 & ( 4, 2, 1, 1 ) &\cr
& -300 & 40 & ( 90, 48, 36, 5, 1 ) & 4 & ( 15, 8, 6, 1 ) &\cr
& -300 & 15 & ( 27, 15, 9, 2, 1 ) & 4 & ( 9, 5, 3, 1 ) &\cr
& -304 & 12 & ( 12, 10, 8, 1, 1 ) & 4 & ( 6, 5, 4, 1 ) &\cr
& -312 & 34 & ( 84, 56, 21, 6, 1 ) & 4 & ( 12, 8, 3, 1 ) &\cr
& -312 & 20 & ( 27, 24, 18, 2, 1 ) & 4 & ( 9, 8, 6, 1 ) &\cr
& -312 & 18 & ( 45, 30, 10, 4, 1 ) & 4 & ( 9, 6, 2, 1 ) &\cr
& -312 & 17 & ( 30, 15, 12, 2, 1 ) & 4 & ( 10, 5, 4, 1 ) &\cr
& -312 & 11 & ( 12, 10, 6, 1, 1 ) & 4 & ( 6, 5, 3, 1 ) &\cr
\noalign{\hrule}
& -312 & 8 & ( 12, 6, 4, 1, 1 ) & 4 & ( 6, 3, 2, 1 ) &\cr
& -312 & 8 & ( 24, 16, 4, 3, 1 ) & 4 & ( 6, 4, 1, 1 ) &\cr
& -312 & 5 & ( 8, 6, 2, 1, 1 ) & 4 & ( 4, 3, 1, 1 ) &\cr
& -336 & 95 & ( 315, 210, 90, 14, 1 ) & 4 & ( 21, 14, 6, 1 ) &\cr
& -336 & 15 & ( 36, 24, 8, 3, 1 ) & 4 & ( 9, 6, 2, 1 ) &\cr
& -348 & 12 & ( 14, 12, 8, 1, 1 ) & 4 & ( 7, 6, 4, 1 ) &\cr
& -360 & 24 & ( 60, 40, 15, 4, 1 ) & 4 & ( 12, 8, 3, 1 ) &\cr
& -360 & 16 & ( 33, 18, 12, 2, 1 ) & 4 & ( 11, 6, 4, 1 ) &\cr
& -360 & 5 & ( 18, 12, 3, 2, 1 ) & 4 & ( 6, 4, 1, 1 ) &\cr
& -368 & 10 & ( 16, 8, 6, 1, 1 ) & 4 & ( 8, 4, 3, 1 ) &\cr
\noalign{\hrule} }}}}$$
\vfill\eject

$${\fivepoint{
\vbox{\offinterlineskip\tabskip=0pt
\halign{\strut\vrule#
&~~$#$~~\hfil
&~~$#$~~\hfil
&~~$#$~~\hfil
&~~$#$~~\hfil
&~~$#$~~\hfil
&\vrule#\cr
\noalign{\hrule}
&Euler~\#& h^{1,1}& \hat w=(w_1,..,w_5) &i &~~~\bar w&\cr
\noalign{\hrule}
& -372 & 80 & ( 273, 182, 78, 12, 1 ) & 4 & ( 21, 14, 6, 1 ) &\cr
& -372 & 36 & ( 105, 70, 28, 6, 1 ) & 4 & ( 15, 10, 4, 1 ) &\cr
& -372 & 8 & ( 14, 8, 4, 1, 1 ) & 4 & ( 7, 4, 2, 1 ) &\cr
& -372 & 4 & ( 10, 6, 2, 1, 1 ) & 4 & ( 5, 3, 1, 1 ) &\cr
& -384 & 25 & ( 60, 32, 24, 3, 1 ) & 4 & ( 15, 8, 6, 1 ) &\cr
& -396 & 32 & ( 90, 60, 24, 5, 1 ) & 4 & ( 15, 10, 4, 1 ) &\cr
& -408 & 18 & ( 48, 32, 12, 3, 1 ) & 4 & ( 12, 8, 3, 1 ) &\cr
& -408 & 10 & ( 27, 18, 6, 2, 1 ) & 4 & ( 9, 6, 2, 1 ) &\cr
& -420 & 10 & ( 18, 10, 6, 1, 1 ) & 4 & ( 9, 5, 3, 1 ) &\cr
& -432 & 59 & ( 210, 140, 60, 9, 1 ) & 4 & ( 21, 14, 6, 1 ) &\cr
\noalign{\hrule}
& -432 & 13 & ( 18, 16, 12, 1, 1 ) & 4 & ( 9, 8, 6, 1 ) &\cr
& -432 & 11 & ( 20, 10, 8, 1, 1 ) & 4 & ( 10, 5, 4, 1 ) &\cr
& -456 & 22 & ( 60, 40, 16, 3, 1 ) & 4 & ( 15, 10, 4, 1 ) &\cr
& -456 & 18 & ( 45, 24, 18, 2, 1 ) & 4 & ( 15, 8, 6, 1 ) &\cr
& -456 & 14 & ( 36, 24, 9, 2, 1 ) & 4 & ( 12, 8, 3, 1 ) &\cr
& -480 & 47 & ( 168, 112, 48, 7, 1 ) & 4 & ( 21, 14, 6, 1 ) &\cr
& -480 & 47 & ( 147, 98, 42, 6, 1 ) & 4 & ( 21, 14, 6, 1 ) &\cr
& -480 & 11 & ( 22, 12, 8, 1, 1 ) & 4 & ( 11, 6, 4, 1 ) &\cr
& -480 & 3 & ( 12, 8, 2, 1, 1 ) & 4 & ( 6, 4, 1, 1 ) &\cr
& -528 & 7 & ( 18, 12, 4, 1, 1 ) & 4 & ( 9, 6, 2, 1 ) &\cr
\noalign{\hrule}
& -552 & 15 & ( 45, 30, 12, 2, 1 ) & 4 & ( 15, 10, 4, 1 ) &\cr
& -564 & 29 & ( 105, 70, 30, 4, 1 ) & 4 & ( 21, 14, 6, 1 ) &\cr
& -612 & 12 & ( 30, 16, 12, 1, 1 ) & 4 & ( 15, 8, 6, 1 ) &\cr
& -624 & 23 & ( 84, 56, 24, 3, 1 ) & 4 & ( 21, 14, 6, 1 ) &\cr
& -624 & 9 & ( 24, 16, 6, 1, 1 ) & 4 & ( 12, 8, 3, 1 ) &\cr
& -720 & 17 & ( 63, 42, 18, 2, 1 ) & 4 & ( 21, 14, 6, 1 ) &\cr
& -732 & 10 & ( 30, 20, 8, 1, 1 ) & 4 & ( 15, 10, 4, 1 ) &\cr
& -960 & 11 & ( 42, 28, 12, 1, 1 ) & 4 & ( 21, 14, 6, 1 ) &\cr
\noalign{\hrule} }}}}$$

\end